%% file: main.tex
\documentclass[12pt]{article}
\usepackage[reqno]{amsmath}
\usepackage{authblk}
\usepackage{amssymb, enumerate,amsthm}
\usepackage{graphicx}
\usepackage{tikz}
\usetikzlibrary{calc}
\usetikzlibrary{decorations.markings}
\usetikzlibrary{decorations.text}
\usepackage{array}
\usepackage{float}
\usepackage{wrapfig}
\usepackage{fullpage}
\usepackage{url}
\usepackage{hyperref}
\usepackage{nicefrac}
\usepackage{blkarray}
\usepackage{dcolumn,booktabs}
\newcolumntype{d}[1]{D{.}{.}{#1}}

\usepackage{nicematrix}
\usepackage{enumitem}
\usepackage[capitalise]{cleveref}
\usepackage{color,colortbl}
\newcommand{\doi}[1]{\href{https://doi.org/#1}{\texttt{doi:#1}}}
\usepackage{kbordermatrix}

\expandafter\ifx\csname pdfpagebox\endcsname\relax\else
\fi
 
\newtheorem{theorem}{Theorem}[section]
\newtheorem*{theorem*}{Theorem}
\newtheorem{lemma}[theorem]{Lemma}

\definecolor{gcolour}{RGB}{230,97,0}
\definecolor{kcolour}{RGB}{93,58,155}

\newcommand\email[1]{\texttt{#1}}

\renewcommand{\leq}{\leqslant}
\renewcommand{\geq}{\geqslant}

\title{Cubic graphs with no eigenvalues in the interval $(-2,0)$ }
\author[1]{Krystal Guo}
\author[2]{Gordon F. Royle}
\affil[1]{Korteweg-de Vries Institute for Mathematics, University of Amsterdam. (\email{k.guo@uva.nl}).}
\affil[2]{Department of Mathematics and Statistics, The University of Western Australia. (\email{gordon.royle@uwa.edu.au}).}
\begin{document}
\maketitle

\begin{abstract}
We give a complete characterisation of the cubic graphs with no eigenvalues in the interval $(-2,0)$. There is one thin infinite family consisting of a single graph on $6n$ vertices for each $n \geqslant 2$, and five ``sporadic'' graphs, namely the $3$-prism $K_3 \mathbin{\square} K_2$, the complete bipartite graph $K_{3,3}$, the Petersen graph, the dodecahedron and Tutte's $8$-cage. The proof starts by observing that if a cubic graph has no eigenvalues in $(-2,0)$ then its local structure around a girth-cycle is very constrained. Then a separate case analysis for each possible girth shows that these constraints can be satisfied only by the known examples. All but one of these case analyses can be completed by hand, but for girth five there are sufficiently many cases that it is necessary to use a computer for the analysis.
\end{abstract}

\input intro

\input easypart

\input hardpart

\input conclusion

\bibliographystyle{plain}
\bibliography{main}

%\printbibliography

\end{document}

%% file: intro.tex
\section{Introduction}

A \emph{spectral gap set} for cubic graphs is an open set $\mathcal{I} \subseteq (-3,3)$ such that there are infinitely many cubic graphs with no eigenvalues in $\mathcal{I}$. Spectral gap sets and the families of graphs realising those gap sets arise naturally in a variety of graph-theoretic contexts. For example, $(2\sqrt{2},3)$ is a spectral gap set realised by cubic expander graphs (see Hoory, Linial and Wigderson \cite{HooLinWig2006} and Chiu \cite{Chi1992}) and spectral gap sets around $0$ are related to the HOMO-LUMO gap of H\"uckel molecular orbital theory and the stability of molecular graphs in chemistry (see Fowler and Pisanski \cite{MR2759780}).

 Koll\'ar and Sarnak \cite{KolSar2021} considered the broader question of identifying further spectral gap sets for cubic graphs, in order to show that every point in $[-3,3)$ is ``gapped'' i.e., belongs to some spectral gap set. As part of this study, they focussed on \emph{maximal} (by inclusion) spectral gap sets and maximal spectral gap \emph{intervals} (i.e., where $\mathcal{I}$ is an open interval). Among numerous other results, they showed that $(-1,1)$ and $(-2,0)$ are maximal spectral gap intervals and that any spectral gap interval for cubic graphs has length at most 2. (We note that Guo and Mohar \cite{GuoMoh2014} had previously shown that $(-1,1)$ is a spectral gap interval, but with a different infinite family of graphs.)

Once a spectral gap set $\mathcal{I}$ has been identified by exhibiting an infinite class of graphs, say $\mathcal{G}$, whose spectrum avoids $\mathcal{I}$, it is natural to ask if the graphs in $\mathcal{G}$ are the \emph{only} such graphs and if not, what are the others? In short, is it possible to give a \emph{complete classification} of the graphs whose spectra avoid $\mathcal{I}$.  For example, the class of line graphs shows that $(-\infty,-2)$ is a spectral gap set, but there are other graphs with least eigenvalue at least $-2$. The complete classification of all such graphs by Cameron, Goethals, Seidel and Shult \cite{CamGoeSeiShu1976} is widely regarded as one of the seminal results of spectral graph theory.

In \cite{guo2024cubicgraphseigenvaluesinterval} we gave a complete classification of the cubic graphs with no eigenvalues in $(-1,1)$, finding two infinite families, along with $14$ sporadic graphs. This classification allowed us to show that $(-1,1)$ is a maximal spectral gap \emph{set} (not just interval). In this paper, we adapt the techniques introduced in \cite{guo2024cubicgraphseigenvaluesinterval} to determine an analogous complete classification of the cubic graphs whose spectrum avoids $(-2,0)$. Although the classification itself is simpler, with just one infinite family and five sporadic graphs, the proof is more complicated and requires more case analysis. Unlike $(-1,1)$, the interval $(-2,0)$ is not a maximal gap set, because the sole infinite family also avoids the intervals  $\left(-3,(-1-\sqrt{17})/2\right)$ and $\left((-1+\sqrt{17})/2,2\right)$ (see  \cref{thm:notmax} and Koll\'ar and Sarnak \cite[Figure~10]{KolSar2021}).

For cubic graphs, the intervals $(-1,1)$ and $(-2,0)$ are the only spectral gap intervals of length 2 with integer endpoints, and there are no further natural spectral gap intervals to analyse. It is conceivable that there may be a spectral gap interval of length 2 with non-integer endpoints, but there is no computational evidence of any such interval.

%(The reader may be as relieved as the authors to learn that $(-1,1)$ and $(-2,0)$ are the only length $2$ spectral gap intervals with integer end-points for which these techniques may work.)

\begin{figure}
\begin{center}
\begin{tikzpicture}
\tikzstyle{vertex}=[fill=blue!25, circle, draw=black, inner sep=0.8mm]
\node [vertex] (v0) at (0,0) {};
\node [vertex] (v1) at (0,-1) {};
\node [vertex] (v2) at (1,0) {};
\node [vertex] (v3) at (1,-1) {};
\node [vertex] (v4) at (2,0) {};
\node [vertex] (v5) at (2,-1) {};
\node [above] at (v0.north) {\small $v_0$};
\node [below] at (v1.south) {\small $v_1$};
\node [above] at (v2.north) {\small $v_2$};
\node [below] at (v3.south) {\small $v_3$};
\node [above] at (v4.north) {\small $v_4$};
\node [below] at (v5.south) {\small $v_5$};
\draw (v0)--(v2);
\draw (v4)--(v3);
\draw (v5)--(v2);
\draw (v1)--(v3);
\draw (v2)--(v4);
\draw (v3)--(v5);
\draw (v0)--(v1);
\draw (v1)--(-0.65,-1);
\draw (v0)--(-0.65,0);
\draw (v4)--(2.65,0);
\draw (v5)--(2.65,-1);
\end{tikzpicture}
\caption{The building block for $X(n)$}
\label{fig:gadget}
\end{center}
\end{figure}
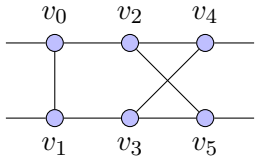

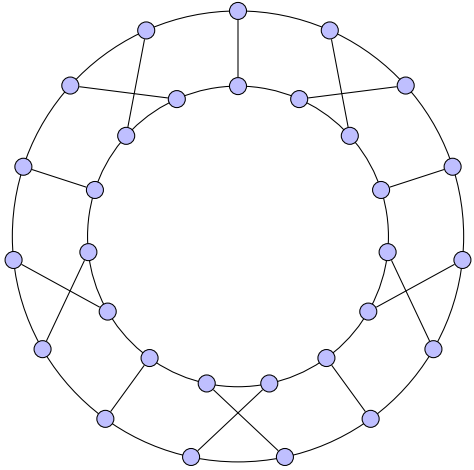
\begin{figure}
\begin{center}
\begin{tikzpicture}
\tikzstyle{vertex}=[fill=blue!25, circle, draw=black, inner sep=0.8mm]
\draw (0,0) circle (2cm);
\draw (0,0) circle (3cm);
\foreach \x in {0,1,...,14} {
  \node [vertex] (v\x) at (90+\x*24:2cm) {};
  \node [vertex] (w\x) at (90+\x*24:3cm) {};
}
\foreach \x in {0,3,6,9,12} {
\draw (v\x)--(w\x);
}
\draw (v1)--(w2);
\draw (v2)--(w1);
\draw (v4)--(w5);
\draw (v5)--(w4);
\draw (v7)--(w8);
\draw (v8)--(w7);
\draw (v10)--(w11);
\draw (v11)--(w10);
\draw (v13)--(w14);
\draw (v14)--(w13);

% \coordinate (a0) at (95:1.8) {};
% \coordinate (b0) at (95:3.2) {};
% \coordinate (c0) at (37:3.2) {};
% \coordinate (d0) at (37:1.8) {};

%\fill [rounded corners = 5pt, black!25!white,opacity=0.5] (a0) -- (b0) [bend left = 30] to (c0) -- (d0) [bend right = 20] to (a0);

%\draw (\ex,\ey) ++(45:.8) arc (45:-45:.8)

\end{tikzpicture}
\caption{$X(5)$ has five gadgets connected in a cyclic fashion}
\label{fig:xk5}

\end{center}
\end{figure}

The infinite family of cubic graphs with no eigenvalues in $(-2,0)$ is constructed as follows: for $n \geqslant 2$, let $X(n)$ be the graph on $6n$ vertices that is constructed by connecting $n$ copies of the $6$-vertex gadget shown in \cref{fig:gadget} in a cyclic fashion as in \cref{fig:xk5}. (The $3$-prism could be viewed as $X(1)$ but it is more convenient to treat the $3$-prism as a sporadic graph and $X(2)$ as the first graph in the infinite family.)

Our main theorem is the following classification:

\begin{theorem}\label{mainresult}
A cubic graph $G$ has no eigenvalues in $(-2,0)$ if and only if $G \cong X(n)$ for some $n \geqslant 2$ or $G$ is the $3$-prism $K_3 \mathbin{\square} K_2$, the complete bipartite graph $K_{3,3}$, the Petersen graph, the dodecahedron, or Tutte's $8$-cage.
\end{theorem}

The remainder of the paper is structured as follows. In \cref{sec:overview} we establish our terminology and notation and give an overview of the overall structure of the proof. In \cref{sec:easypart} we give the ``easy part'' of the proof by showing that the graphs in our list have no eigenvalues in $(-2,0)$. Finally, we prove that these are the \emph{only} cubic graphs with no eigenvalues in $(-2,0)$ in a sequence of subsections, one for each possible girth. Finally we end with some concluding remarks.

\input outline

%\input voltages

%% file: outline.tex
\section{Overview}\label{sec:overview}

Throughout the remainder of this paper $G$ will denote a cubic graph and $A(G)$ (or just $A$) its adjacency matrix.

We split the proof into two parts, first the relatively straightforward task of showing that the graphs listed in \cref{mainresult} have no eigenvalues in $(-2,0)$, and then the much harder task of showing that these are the \emph{only} cubic graphs with no eigenvalues in $(-2,0)$.

For this second task, the matrix $M(G) = A(G) (A(G) + 2 I)$ (or just $M$) will play a major role in the arguments that follow.  If $v$ is an eigenvector of $A$ such that $Av = \lambda v$, then $Mv = \lambda (\lambda+2) v$ and so $v$ is also an eigenvector of $M$. If $A$ has no eigenvalues in $(-2,0)$, then every eigenvalue of $M$ is non-negative and so $M$ is positive semidefinite. 

% \krystalsays{Any principal submatrix of a positive semidefinite matrix is also positive semidefinite. Thus the existence of a principal submatrix with at least one negative eigenvalue is sufficient to show that $M$ is not positive semidefinite. If a principal submatrix has negative determinant, then it has an odd number of negative eigenvalues and is thus not positive semidefinite. We may want this statement instead because it could be, in some case, it will be convenient to exhibit a principal submatrix of $M$ which is not positive semidefinite, but has determinant $0$ (because it has 2 negative eigenvalues).}

We recall that a matrix is positive semidefinite if and only if all of its principal minors are non-negative. Therefore the existence of \emph{any} principal submatrix of $M$ with negative determinant suffices to show that $M$ is \emph{not} positive semidefinite. We repeatedly use this observation by showing that certain (smallish) configurations in a cubic graph would determine a negative principal minor of $M$ and therefore cannot occur in a graph with no eigenvalues in $(-2,0)$. The connection between subgraphs in $G$ and principal minors of $M$ arises from the following combinatorial interpretation of the entries of $M$: 
\[
M_{uv} = 
\begin{cases}
3, & u = v;\\
|N(u) \cap N(v)|, & u \not\sim v;\\
2 + |N(u) \cap N(v)|, & u \sim v.\\
\end{cases}
\]
where $N(u)$ denotes the open neighbourhood of $u$.

\newcommand\cor[1]{\mathrm{Cor}(#1)}

The final stage of the characterisation is to show that the only cubic graphs that avoid these forbidden configurations are the ones listed in \cref{mainresult}. 

Finally, we note that although our lemmas are phrased in a \emph{qualitative} fashion, they can be rewritten in a more precise \emph{quantitative} fashion. Rather than simply asserting that a certain configuration is forbidden in a cubic graph with no eigenvalues in $(-2,0)$, we can say that any cubic graph \emph{with} a certain configuration has an eigenvalue in a specific smaller subinterval of $(-2,0)$. 

% In the various lemmas that follow, we identify configurations that are incompatible with graphs with no 
% We note that the  a negative principal minor of $M$ 
% In fact, we can show if $G$ is not one of the graphs in our classification, $G$ has a root in a smaller subinterval of $(-2,0)$ centered at $-1$. 
% %For example, if $G$ is a cubic graph of girth $3$ and  $G$ is not isomorphic to the $3$-prism, then $G$ has an eigenvalue in $[-1.893358, -0.106642]$. 
% To establish these results, we need the following. 

\begin{lemma}\label{lem:mtau}
If $M'$ is a principal submatrix of $M =A(A+2I)$ with least eigenvalue $\tau'$, then $A$ has an eigenvalue in the closed interval 
$[	-1 - \sqrt{1+ \tau'}, 	-1 + \sqrt{1+ \tau'}  ]. $
%Further, if $M'$ has least eigenvalue $-1$, then $A$ has eigenvalue $-1$. 
\end{lemma}

\begin{proof}
If $v$ is an eigenvector for $A$ such that $Av = \lambda v$, then $Mv = \lambda(\lambda+2)v$ and so $v$ is an eigenvector for $M$ with eigenvalue $\lambda(\lambda+2)$. Conversely, if $v$ is an eigenvector for $M$ with eigenvalue $\tau$, then it is an eigenvector for $A$ with eigenvalue $\lambda = -1 \pm \sqrt{1+\tau}$ (one of the two solutions to $\lambda(\lambda+2)=\tau$).

If $M$ has least eigenvalue $\tau$, then
\[
-1 \leqslant \tau \leqslant \tau'
\]
where the first inequality holds because the minimum value of $\lambda(\lambda+2)$ is $-1$ achieved only when $\lambda = -1$, and the second holds by interlacing. It follows that $0 \leq 1+\tau \leq 1+\tau'$ and so 
\[
-1-\sqrt{1+\tau'} \leq -1 -\sqrt{1+\tau}  \leq -1 + \sqrt{1+\tau} \leq 1+\sqrt{1+\tau'}.
\]
and so the given interval contains both of the values $-1 \pm \sqrt{1+\tau}$.
\end{proof}

We illustrate this concept with a specific example in \cref{lem:quant}.

% The matrices $A$ and $M$ have the same eigenvectors and each eigenvalue $\lambda$ of $A$ corre
% Suppose that $\v$ is an eigenvector of $M$ for $\tau$ in this eigenbasis. Let $\lambda$ be the corresponding eigenvalue of $A$; $A\Zv = \lambda \Zv$. We see that 
% \[
% \lambda(\lambda+2) = \tau.
% \]
% We have then $\lambda^2 + 2\lambda - \tau=0$ and then 
% \[
% \lambda = \frac{-2 \pm \sqrt{4+4\tau}}{2} = -1 \pm \sqrt{1+\tau}. 
% \]
% Since all eigenvalues of $A$ are real, we see that $\tau \geq -1$. 
% Thus $A$ has a root in the set $\{-1 -\sqrt{1+\tau}, -1 + \sqrt{1+\tau}\}$. 

% If $M'$ is a principal submatrix of $M$ with least eigenvalue $\tau'$, we have that $\tau' \geq \tau$ by interlacing.
% \[
% \begin{split}
% 	\tau' &\geq \tau \geq -1  \\
% 	1+ \tau' &\geq 1+\tau \geq 0  \\
% 	\sqrt{1+ \tau'} &\geq \sqrt{ 1+\tau} \geq 0  \\
% \end{split}
% \]
% and thus 
% \[
% 	-1 + \sqrt{1+ \tau'} \geq -1 + \sqrt{ 1+\tau}
% \]
% and 
% \[
% 	-1 - \sqrt{1+ \tau'} \leq -1 - \sqrt{ 1+\tau}. 
% \]
% Thus 
% \[
% \{-1 -\sqrt{1+\tau}, -1 + \sqrt{1+\tau}\} \in [	-1 - \sqrt{1+ \tau'}, 	-1 + \sqrt{1+ \tau'}  ]
% \]
% and the statement follows. If $\tau' = -1$, then $\tau = -1$ since $A$ is real-rooted and has $-1$ as a root. 
% \end{proof}

%% file: easypart.tex
\section{The graphs}\label{sec:easypart}

In this section we describe the graphs mentioned in \cref{mainresult} and verify that they have no eigenvalues in $(-2,0)$.

First we consider the five sporadic graphs, which have the following spectra:
\begin{align*}
\mathrm{sp}(K_3 \mathbin{\square} K_2) &= \{3, 1, 0^{(2)}, -2^{2}\},\\
\mathrm{sp}({K_{3,3}})&= \{3, 0^{(4)}, -3\},\\
\mathrm{sp}(P)&= \{3, 1^{(5)}, -2^{(4)} \}, \\
\mathrm{sp}(D)&= \{3, \sqrt{5}^{(3)},  0^{(4)}, -2^{(4)}, -\sqrt{5}^{(3)}, 3\},\\
\mathrm{sp}(T)&= \{3, 2^{(9)}, 0^{(10)}, -2^{(9)}, -3\}.
\end{align*}
By inspection these contain no eigenvalues in $(-2,0)$.

Next we show that the family $X(n)$ has no eigenvalues in $(-2,0)$.  We use a well-known result from matrix theory.

\begin{theorem}\label{thm:blockcirc}
If $M$ is a block-circulant matrix with $n$ square blocks $A_0$, $A_1$, $\ldots$, $A_{n-1}$, then 
\[
\det M = \prod_{\omega^n = 1} \det \left( A_0 + \omega A_1 + \cdots + \omega^{n-1} A_{n-1} \right),
\]
where $\omega$ ranges over the complex $n$-th roots of unity. \qed
\end{theorem}

The adjacency matrix of $X(n)$ is an $n \times n$ block-circulant matrix of the form
\[
\renewcommand{\arraystretch}{1.6}
A(X(n)) = \left[\begin{array}{ccccc}
A&B&&&B^T\\
B^T&A&B&&\\
&\ddots&\ddots&\ddots&\\
&&B^T&A&B\\
B&&&B^T&A\\
\end{array}
\right]
\]
with $6 \times 6$ blocks, of which the only ones that are non-zero are 
\[
A =
\kbordermatrix{
& v_0 & v_1 & v_2 & v_3 & v_4 & v_5 \\
v_0 & 0 & 1 & 1 & 0 & 0 & 0\\
v_1 & 1 & 0 & 0 & 1 & 0 & 0\\
v_2 & 1 & 0 & 0 & 0 & 1 & 1\\
v_3 & 0 & 1& 0 & 0 & 1 & 1\\
v_4 & 0 & 0 & 1 & 1 & 0 & 0\\
v_5 & 0 & 0 & 1 & 1 & 0 & 0\\
},
\qquad 
B = \left[
\begin{array}{rrrrrr}
0&0&0&0&0&0\\
0&0&0&0&0&0\\
0&0&0&0&0&0\\
0&0&0&0&0&0\\
1&0&0&0&0&0\\
0&1&0&0&0&0
\end{array}
\right],
\]
and $B^T$.
Here $A$ accounts for the edges inside each $6$-vertex gadget, while $B$ and $B^T$ account for the two edges connecting each gadget to the previous and subsequent gadgets in the cyclic order.

\begin{theorem}
The characteristic polynomial of $X(n)$ is given by
\[
x^n (x-1)^n (x+2)^n 
\prod_{j=0}^{j=n-1}
\left(x^3-x^2-6 x+4 \left( 1 - \cos \frac{2\pi j}{n} \right)\right).
\]
\end{theorem}

\crefname{expr}{Expression}{Expression}

\begin{proof}
The characteristic polynomial of $X(n)$ is $\det (xI - A(X(n)))$ which is also a block-circulant matrix. Using the formula from \cref{thm:blockcirc} and noting that $\omega^{n-1} = \omega^{-1}$ it follows that $\varphi(X(n))$ is the product of the $n$ determinants
\begin{equation}\label{eq:det}
\left|
\begin{array}{rrrrrr}
 x & -1 & -1 & 0 & -\frac{1}{\omega} & 0 \\
 -1 & x & 0 & -1 & 0 & -\frac{1}{\omega} \\
 -1 & 0 & x & 0 & -1 & -1 \\
 0 & -1 & 0 & x & -1 & -1 \\
 -\omega & 0 & -1 & -1 & x & 0 \\
 0 & -\omega & -1 & -1 & 0 & x \\
\end{array}
\right|
\end{equation}

as $\omega$ ranges over the complex $n$'th roots of unity. If we substitute  $\omega = \cos \frac{2\pi j}{n} + i \sin \frac{2\pi j}{n}$ (where $0 \leq j < n$), the determinant given in \eqref{eq:det} evaluates to\
\[
x (x-1)(x+2) \left(x^3-x^2-6 x+4 \left( 1 - \cos \frac{2\pi j}{n} \right)\right),
\]
and so the result holds.
\end{proof}

\begin{theorem}
For $n \geqslant 2$, the graph $X(n)$ has no eigenvalues in the three intervals $(-2,0)$, $\left(-3,(-1-\sqrt{17})/2\right)$ and $\left((-1+\sqrt{17})/2,2\right)$.
\end{theorem}

\begin{proof}\label{thm:notmax}
The eigenvalues of $X(n)$ other than $-2$, $0$ and $1$ are the roots of cubic polynomials of the form $x^3-x^2-6x+\alpha$ where $\alpha \in [0,8]$. It is easy to verify that $x^3-x^2-6x+8 < 0$ when $x \in \left(-3,(-1-\sqrt{17})/2\right)$ or when $x \in \left((-1+\sqrt{17})/2,2\right)$, and that $x^3-x^2-6x > 0$ for $x \in (-2,0)$.
\end{proof}

%% file: hardpart.tex
\section{The characterisation}\label{sec:hardpart}

In this long section we break down the characterisation of cubic graphs with no eigenvalues in $(-2,0)$ into seven subsections, one each for girths $3$, $4$, $5$, $6$, $7$, $8$ and a final subsection for girth greater than $8$. 

As previously, $G$ denotes a cubic graph with adjacency matrix $A$ and $M = A(A+2I)$. 

In general, the arguments show that certain subgraphs cannot occur in $G$ by showing that it would then possible to find a subset $T \subseteq V(G)$ such that $M_{TT}$ has negative determinant. Because variants of this argument are used repeatedly in what follows, we emphasize that the entries of $M_{TT}$ are not determined solely by the subgraph induced by $T$, but also by any common neighbours of pairs of vertices in $T$. 

We rely heavily on analysing the structure of the graph close to a girth cycle, so we establish some precise terminology to describe this situation: If $C = (u_0,u_1,\ldots,u_{g-1})$ is a cycle of length $g$ in a cubic graph $G$ of girth $g$, then each vertex $u_i$ has a unique third neighbour $v_i$ that is not on $C$. The \emph{corona} of $C$ is the subgraph $\cor{C}$ with
\begin{align*}
V(\cor{C}) &= \{u_0,u_1,\ldots,u_{g-1},v_0,v_1,\ldots,v_{g-1}\},\\
E(\cor{C}) &= E(C) \cup \{u_iv_i : 0 \leq i < g\}.
\end{align*}
If the vertices $V = \{v_0,v_1,\ldots,v_{g-1}\}$ happen to be distinct and pairwise non-adjacent then $\cor{C}$ is an induced subgraph isomorphic to the \emph{corona product} $C_g \circ K_1$ (see Frucht and Harary \cite{MR281659} for the definition and \cref{fig:corona5} for an example), but in general $\cor{C}$ is not an induced subgraph of $G$.

\input girth3

\input girth4NEW

\input girth5X

\subsection{Girth six}

This subsection disposes of the case where the girth is $6$.

\begin{lemma}
A cubic graph $G$ of girth $6$ has an eigenvalue in the interval $(-2,0)$.
\end{lemma}
\begin{proof}

Let $C$ be a $6$-cycle and consider its corona $\cor{C}$ labelled as in \cref{fig:girth6}. As $G$ has girth $6$, the vertices 
$V = \{v_0,v_1,\ldots,v_5\}$ are distinct and each vertex $v_i \in V$ may be adjacent to $v_{i \pm 3}$ but cannot be adjacent to any other vertex in $V$. If we take $T = \{u_0, u_1, u_3, u_4, v_0\}$ (as highlighted in the diagram), then 
\[
M_{TT} = 
\kbordermatrix{
& u_0 & u_1 & u_3 & u_4 & v_0 \\
u_0 & 3 & 2 & 0 & 1 & 2 \\
u_1 & 2 & 3 & 1 & 0 & 1 \\
u_3 &0 & 1 & 3 & 2 & \alpha \\
u_4 &1 & 0 & 2 & 3 & 0 \\
v_0 & 2 & 1 & \alpha & 0 & 3 \\
}
\]
% \[
% M_{TT} = 
% \left[
% \begin{array}{ccccc}
% 3 & 2 & 0 & 1 & 2 \\
% 2 & 3 & 1 & 0 & 1 \\
% 0 & 1 & 3 & 2 & \alpha \\
% 1 & 0 & 2 & 3 & 0 \\
% 2 & 1 & \alpha & 0 & 3 \\
% \end{array}
% \right]
% \]
where $\alpha$ is either $0$ or $1$ depending on whether $v_0$ is adjacent to $v_3$ or not. 

This matrix has determinant $-48$ if $\alpha=0$ and $-12$ if $\alpha=1$, and so in either case $M$ is not positive semidefinite, and we conclude that $G$ has an eigenvalue in $(-2,0)$.
\end{proof}

\begin{figure}
\begin{center}
\begin{tikzpicture}
\tikzstyle{vertex}=[fill=blue!25, circle, draw=black, inner sep=0.8mm]
\tikzstyle{rvertex}=[fill=red, circle, draw=black, inner sep=0.8mm]
\foreach \x in {0,1,...,5} {
\node [vertex,label={[label distance=-0.85cm]60*\x:$v_\x$}] (w\x) at (60*\x:1.75cm) {};
\node [vertex,label={[label distance=-1mm]60*\x:$u_\x$}] (v\x) at (60*\x:2.75cm) {};
\draw (v\x)--(w\x);
}
\draw (v0)--(v1)--(v2)--(v3)--(v4)--(v5)--(v0);
\node [rvertex] at (w0) {};
\node [rvertex] at (v0) {};
\node [rvertex] at (v1) {};
\node [rvertex] at (v3) {};
\node [rvertex] at (v4) {};
\draw [ultra thick, red] (w0)--(v0)--(v1);
\draw [ultra thick, red] (v3)--(v4);
\end{tikzpicture}
\end{center}
\caption{Corona of a $6$-cycle in a graph of girth $6$.}
\label{fig:girth6}
\end{figure}

\subsection{Girth seven}

This subsection disposes of the case where the girth is $7$.

\begin{lemma}
A cubic graph $G$ of girth $7$ has an eigenvalue in the interval $(-2,0)$.
\end{lemma}
\begin{proof}
Let $C$ be a $7$-cycle and consider its corona $\cor{C}$ labelled as in \cref{fig:girth7}. As $G$ has girth $7$, the vertices $V = \{v_0,v_1,\ldots,v_6\}$ are distinct and pairwise non-adjacent. Each vertex $v_i \in V$ may have a common neighbour with $v_{i \pm 3}$ but not with any of the other vertices of $V$. Therefore, if we take $T = \{u_0, u_1,u_3,u_5,v_0,v_5\}$ then

\[
M_{TT} = 
\kbordermatrix{
&u_0& u_1&u_3&u_5&v_0&v_5\\
u_0&3&2&0&1&2&0\\
u_1&2&3&1&0&1&0\\
u_3&0&1&3&1&0&0\\
u_5&1&0&1&3&0&2\\
v_0&2&1&0&0&3&0\\
v_5&0&0&0&2&0&3}
\]
which has determinant $-16$. 
% \[
% M_{TT} = 
% \left[
% \begin{array}{cccccc}
% 3&2&0&1&2&0\\
% 2&3&1&0&1&0\\
% 0&1&3&1&0&0\\
% 1&0&1&3&0&2\\
% 2&1&0&0&3&0\\
% 0&0&0&2&0&3\\
% \end{array}
% \right]
% \]

\end{proof}

\begin{figure}
\begin{center}
\begin{tikzpicture}
\tikzstyle{vertex}=[fill=blue!25, circle, draw=black, inner sep=0.8mm]
\tikzstyle{rvertex}=[fill=red, circle, draw=black, inner sep=0.8mm]
\foreach \x in {0,1,...,6} {
\node [vertex](v\x) at (90+51.4*\x:1.75cm) {};
\node at (90+51.4*\x:1.35cm) {$v_{\x}$};

\node [vertex](u\x) at (90+51.4*\x:2.75cm) {};
\node at (90+51.4*\x:3.2cm) {$u_{\x}$};
\draw (u\x)--(v\x);
}
\draw (u0)--(u1)--(u2)--(u3)--(u4)--(u5)--(u6)--(u0);
\node [rvertex] at (u0) {};
\node [rvertex] at (u1) {};
\node [rvertex] at (u3) {};
\node [rvertex] at (u5) {};
\node [rvertex] at (v0) {};
\node [rvertex] at (v5) {};
\draw [ultra thick, red] (v0)--(u0)--(u1);
\draw [ultra thick, red] (u5)--(v5);
\end{tikzpicture}
\end{center}
\caption{Corona of a $7$-cycle in a graph of girth $7$.}
\label{fig:girth7}
\end{figure}

\subsection{Girth eight}

In this subsection we will prove that the only cubic graph of girth $8$ with no eigenvalues in $(-2,0)$ is Tutte's $8$-cage, which is the unique cubic graph on $30$ vertices of girth $8$. 

We start with the corona of an $8$-cycle labelled as in \cref{fig:corona8}.  As $G$ has girth $8$, the vertices $V = \{v_0, v_1, \ldots, v_7\}$ are distinct and pairwise non-adjacent. Each vertex $v_i$ may have a common neighbour with $v_{i+4}$ but not with any of the other vertices of $V$.

\begin{figure}
\begin{center}
\begin{tikzpicture}
\tikzstyle{vertex}=[fill=blue!25, circle, draw=black, inner sep=0.8mm]
\tikzstyle{rvertex}=[fill=red, circle, draw=black, inner sep=0.8mm]
\foreach \x in {0,1,...,7} {
\node [vertex,label={[label distance=-1cm]90+45*\x:$v_\x$}] (v\x) at (90+45*\x:1.75cm) {};
\node [vertex,label={[label distance=0cm]90+45*\x:$u_\x$}] (u\x) at (90+45*\x:2.75cm) {};
\draw (u\x)--(v\x);
}
\draw (u0)--(u1)--(u2)--(u3)--(u4)--(u5)--(u6)--(u7)--(u0);

\pgftransformxshift{9cm}
\foreach \x in {0,1,...,7} {
\node [vertex,label={[label distance=-1cm]45*\x:$v_\x$}] (v\x) at (45*\x:1.75cm) {};
\node [vertex,label={[label distance=0cm]45*\x:$u_\x$}] (u\x) at (45*\x:2.75cm) {};
\draw (u\x)--(v\x);
}
\draw (u0)--(u1)--(u2)--(u3)--(u4)--(u5)--(u6)--(u7)--(u0);

\node [rvertex] at (u0) {};
\node [rvertex] at (v0) {};
\node [rvertex] at (u2) {};
\node [rvertex] at (v2) {};
\node [rvertex] at (u4) {};
\node [rvertex] at (v4) {};
\node [rvertex] at (u6) {};
\draw [ultra thick, red] (u0)--(v0);
\draw [ultra thick, red] (u2)--(v2);
\draw [ultra thick, red] (u4)--(v4);
\end{tikzpicture}
\caption{Corona of an $8$-cycle in a graph of girth $8$.}
\label{fig:corona8}
\end{center}
\end{figure}

If we take $T = \{u_0, v_0, u_2, v_2, u_4, v_4, u_6\}$, then the   submatrix $M_{TT}$ is given by
\[
M_{TT} = 
\kbordermatrix{
& u_0& v_0& u_2& v_2& u_4& v_4& u_6\\
u_0&3&2&1&0&0&0&1\\
v_0&2&3&0&0&0&\alpha&0\\
u_2&1&0&3&2&1&0&0\\
v_2&0&0&2&3&0&0&0\\
u_4&0&0&1&0&3&2&1\\
v_4&0&\alpha&0&0&2&3&0\\
u_6&1&0&0&0&1&0&3
}
\]
where $\alpha = 1$ if $v_0$ and $v_4$ have a common neighbour and $0$ otherwise. If $\alpha = 0$ then $\det M_{TT} = -45$ and so $G$ has an eigenvalue in $(-2,0)$. 

Therefore if $G$ has no eigenvalues in $(-2,0)$, it must be the case that $v_0$ and $v_4$ have a common neighbour, denote this $v_{04}$.  By symmetry, it follows that $\{v_1, v_5\}$ must have a common neighbour $v_{15}$ and similarly we get $v_{26}$ and $v_{37}$. As $G$ has girth $8$, the vertices $\{v_{04},v_{15},v_{26},v_{37}\}$ are distinct, and thus $G$ contains the $20$-vertex configuration shown in \cref{fig:20vertex}. 

As we will repeatedly be using it, we emphasize that this argument implies that if $C$ is \emph{any} $8$-cycle in $G$ and $s$ and $t$ are antipodal vertices on $C$, then their neighbours $s' \sim s$ and $t' \sim t$ off $C$ must themselves share a common neighbour. 

\begin{figure}
\begin{center}
\begin{tikzpicture}
\tikzstyle{vertex}=[fill=blue!25, circle, draw=black, inner sep=0.8mm]
\tikzstyle{rvertex}=[fill=red, circle, draw=black, inner sep=0.6mm]
\foreach \x in {0,1,...,7} {
\node [vertex] (v\x) at (45*\x:2.5cm) {};
\node at (45*\x:2.1cm) {\small $v_\x$};
\node [vertex] (u\x) at (45*\x:4cm) {};
\node at (45*\x:4.4cm) {\small $u_\x$};
\draw (u\x)--(v\x);
}
\draw (u0)--(u1)--(u2)--(u3)--(u4)--(u5)--(u6)--(u7)--(u0);

\node [vertex,label=above:{\small $v_{04}$}] (v04) at (0,1) {};
\node [vertex,label=below right:{\small $v_{26}$}] (v26) at (1,0) {};
\node [vertex,label=below:{\small $v_{15}$}] (v15) at (0,-1) {};
\node [vertex,label=below:{\small $v_{37}$}] (v37) at (-1,0) {};
\draw (v0)--(v04)--(v4);
\draw (v2)--(v26)--(v6);
\draw (v3)--(v37)--(v7);
\draw (v1)--(v15)--(v5);

% \foreach \x in {0,1,...,7} {
%   \node [vertex,label={[label distance=-0.9cm]45*\x:$x_\x$}] (x\x) at (22.5+45*\x:3.25cm) {};
%   \draw (w\x) -- (x\x);
% }

% \draw [line width=6pt, yellow, opacity=0.5] (v0.center)--(v1.center)--(v2.center)--(v3.center)--(w3.center)--(w37.center)--(w7.center)--(v7.center)--(v0.center);

% \draw [line width=6pt, red!50, opacity=0.5] (v0.center)--(w0.center)--(w04.center)--(w4.center)--(v4.center)--(v5.center)--(v6.center)--(v7.center)--(v0.center);

%\draw (x1) [bend right = 80] to (x4) [bend right = 80] to (x7) [bend right = 80] to (x2) [bend right = 80] to (x5) [bend right = 80] to (x0) [bend right = 80] to (x3) [bend right = 80] to (x6) [bend right = 80] to (x1);

% \draw (x1) .. controls (112.5:4.75cm) and (157.5:4.75cm) ..(x4);
% \draw (x2) .. controls (157.5:4.75cm) and (202.5:4.75cm) ..(x5);
% \draw (x3) .. controls (202.5:4.75cm) and (247.5:4.75cm) ..(x6);
% \draw (x4) .. controls (247.5:4.75cm) and (292.5:4.75cm) ..(x7);
% \draw (x5) .. controls (292.5:4.75cm) and (337.5:4.75cm) ..(x0);
% \draw (x6) .. controls (337.5:4.75cm) and (22.5:4.75cm) ..(x1);
% \draw (x7) .. controls (22.5:4.75cm) and (67.5:4.75cm) ..(x2);
% \draw (x0) .. controls (67.5:4.75cm) and (112.5:4.75cm) ..(x3);

\end{tikzpicture}
\caption{Structure induced by an $8$-cycle}
\label{fig:20vertex}
\end{center}
\end{figure}

The remainder of the argument is a case-analysis starting from the $20$-vertex configuration and showing that it can only be completed to Tutte's 8-cage.  Our next step is to consider the third neighbours, say $\{v_0,v_1,\ldots,v_7\}$ of $\{u_0, u_1, \ldots, u_7\}$. These new vertices must be distinct and also not equal to any vertex in $\{v_{04},v_{15},v_{26},v_{37}\}$ because otherwise $G$ would contain a cycle of length less than $8$.

\begin{figure}
\begin{center}
\begin{tikzpicture}[scale=0.85]
\tikzstyle{vertex}=[fill=blue!25, circle, draw=black, inner sep=0.6mm]
\tikzstyle{rvertex}=[fill=red, circle, draw=black, inner sep=0.8mm]
\foreach \x in {0,1,...,7} {
\node [vertex] (v\x) at (45*\x:2.5cm) {};
\node at (45*\x:2.1cm) {\small $v_\x$};
\node [vertex] (u\x) at (45*\x:4cm) {};
\node at (45*\x:4.4cm) {\small $u_\x$};
\draw (u\x)--(v\x);
\node [vertex] (w\x) at (15+45*\x:2.9cm) {};
\node at (15+45*\x:3.3cm) {\small $w_\x$};
\draw (v\x)--(w\x);

}
\draw (u0)--(u1)--(u2)--(u3)--(u4)--(u5)--(u6)--(u7)--(u0);

\node [vertex,label=above:{\small $v_{04}$}] (v04) at (0,1) {};
\node [vertex,label=below right:{\small $v_{26}$}] (v26) at (1,0) {};
\node [vertex,label=below:{\small $v_{15}$}] (v15) at (0,-1) {};
\node [vertex,label=below:{\small $v_{37}$}] (v37) at (-1,0) {};
\draw (v0)--(v04)--(v4);
\draw (v2)--(v26)--(v6);
\draw (v3)--(v37)--(v7);
\draw (v1)--(v15)--(v5);

\draw [line width=5pt, opacity=0.5, red, rounded corners] (u0.center)--(v0.center)--(v04.center)--(v4.center)--(u4.center)--(u5.center)--(u6.center)--(u7.center)--cycle;

\node [rvertex] at (v04) {};
\node [rvertex] at (v4) {};
\node [rvertex] at (u4) {};
\node [rvertex] at (u5) {};
\node [rvertex] at (u6) {};
\node [rvertex] at (u7) {};
\node [rvertex] at (u0) {};
\node [rvertex] at (v0) {};

\end{tikzpicture}
\caption{A second $8$-cycle in the graph}
\label{fig:partialtutte0}
\end{center}
\end{figure}

% We know that if $C$ is any $8$-cycle in $G$, then each vertex has a unique neighbour off $C$ and the unique neighbours of opposite vertices of $C$ must themselves have a common neighbour. 

\begin{figure}
\begin{center}
\begin{tikzpicture}[scale=0.85]
\tikzstyle{vertex}=[fill=blue!25, circle, draw=black, inner sep=0.6mm]
\tikzstyle{rvertex}=[fill=red, circle, draw=black, inner sep=0.8mm]
\foreach \x in {0,1,...,7} {
\node [vertex] (v\x) at (45*\x:2.5cm) {};
\node at (45*\x:2.1cm) {\small $v_\x$};
\node [vertex] (u\x) at (45*\x:4cm) {};
\node at (45*\x:4.4cm) {\small $u_\x$};
\draw (u\x)--(v\x);

\node [vertex] (w\x) at (15+45*\x:2.9cm) {};
\node at (15+45*\x:3.3cm) {\small $w_\x$};

\coordinate  (wp\x) at (15+45*\x:5.75cm) {};

\draw (v\x)--(w\x);
% \node (wp\x) at (10+45*\x:3.5cm) {};
% \node (wm\x) at (20+45*\x:3.5cm) {};
% \draw (w\x)--(wp\x);
% \draw (w\x)--(wm\x);

}
\draw (u0)--(u1)--(u2)--(u3)--(u4)--(u5)--(u6)--(u7)--(u0);

\draw (w0) .. controls (wp1) and (wp2) .. (w3);
\draw (w1) .. controls (wp2) and (wp3) .. (w4);
\draw (w2) .. controls (wp3) and (wp4) .. (w5);
\draw (w3) .. controls (wp4) and (wp5) .. (w6);
\draw (w4) .. controls (wp5) and (wp6) .. (w7);
\draw (w5) .. controls (wp6) and (wp7) .. (w0);
\draw (w6) .. controls (wp7) and (wp0) .. (w1);
\draw (w7) .. controls (wp0) and (wp1) .. (w2);

\node [vertex,label=above:{\small $v_{04}$}] (v04) at (0,1) {};
\node [vertex,label=below right:{\small $v_{26}$}] (v26) at (1,0) {};
\node [vertex,label=below:{\small $v_{15}$}] (v15) at (0,-1) {};
\node [vertex,label=below:{\small $v_{37}$}] (v37) at (-1,0) {};
\draw (v0)--(v04)--(v4);
\draw (v2)--(v26)--(v6);
\draw (v3)--(v37)--(v7);
\draw (v1)--(v15)--(v5);

% \node [vertex,label=above right:{\small $a$}] (a) at (0.25,0.25) {};
% \node [vertex,label=below left:{\small $b$}] (b) at (-0.125,-0.125){};

% \draw (v04)--(a)--(v26);
% \draw (v15)--(b)--(v37);

\draw [line width=5pt, opacity=0.5, red, rounded corners] (u0.center)--(v0.center)--(v04.center)--(v4.center)--(u4.center)--(u5.center)--(u6.center)--(u7.center)--cycle;
%\draw (w0) -- (w3) -- (w6)--(w1)--(w4)--(w7)--(w2)--(w5)--(w0);

\node [rvertex] at (v04) {};
\node [rvertex] at (v4) {};
\node [rvertex] at (u4) {};
\node [rvertex] at (u5) {};
\node [rvertex] at (u6) {};
\node [rvertex] at (u7) {};
\node [rvertex] at (u0) {};
\node [rvertex] at (v0) {};

\end{tikzpicture}
\caption{A $28$-vertex graph with four vertices of degree two}
\label{fig:partialtutte1}
\end{center}
\end{figure}

Now we will use the argument about antipodal vertices on an $8$-cycle, in particular we consider the $8$-cycle
\[
D = (u_0,v_0,v_{04},v_4,u_4,u_5,u_6,u_7)
\]
as highlighted in red in \cref{fig:partialtutte0}.

On the cycle $D$ the vertex $v_4$ is antipodal to $u_7$ and these vertices have ``third neighbours'' $w_4$ and $v_7$. Therefore $w_4$ and $v_7$ must have a common neighbour and the only possibility for this common neighbour is $w_7$, so $w_4 \sim w_7$. The graph has $8$-fold rotational symmetry permuting the indices modulo $8$, and so every vertex $w_i$ is adjacent to $w_{i\pm 3}$.  By this stage we have constructed a $28$-vertex graph, depicted in \cref{fig:partialtutte1} where every vertex has degree $3$ other than $\{v_{04},v_{15},v_{26},v_{37}\}$. 

The girth constraint means that $v_{04}$ is not adjacent to any of the vertices $\{v_{15},v_{26},v_{37}\}$ and so its third neighbour is a new vertex, denote this $a$. Returning to the red $8$-cycle $D$ (redrawn in \cref{fig:partialtutte1} we see that $v_{04}$ and $u_6$ are antipodal on $D$ and therefore their third neighbours, namely $a$ and $v_6$, must share a common neighbour. Since the only neighbour of $v_6$ that does not already have degree $3$ is $v_{26}$ it must be case that $a \sim v_{26}$. By symmetry there is another new vertex $b$ with neighbours $v_{15}$ and $v_{37}$.  So now we have constructed a 30-vertex graph that is Tutte's $8$-cage minus a single edge, as depicted in \cref{fig:partialtutte2}. To complete the argument we must show that $a \sim b$. 

\begin{figure}
\begin{center}
\begin{tikzpicture}[scale=0.85]
\tikzstyle{vertex}=[fill=blue!25, circle, draw=black, inner sep=0.6mm]
\tikzstyle{rvertex}=[fill=red, circle, draw=black, inner sep=0.8mm]
\foreach \x in {0,1,...,7} {
\node [vertex] (v\x) at (45*\x:2.5cm) {};
\node at (45*\x:2.1cm) {\small $v_\x$};
\node [vertex] (u\x) at (45*\x:4cm) {};
\node at (45*\x:4.4cm) {\small $u_\x$};
\draw (u\x)--(v\x);

\node [vertex] (w\x) at (15+45*\x:2.9cm) {};
\node at (15+45*\x:3.3cm) {\small $w_\x$};

\coordinate  (wp\x) at (15+45*\x:5.75cm) {};

\draw (v\x)--(w\x);
% \node (wp\x) at (10+45*\x:3.5cm) {};
% \node (wm\x) at (20+45*\x:3.5cm) {};
% \draw (w\x)--(wp\x);
% \draw (w\x)--(wm\x);

}
\draw (u0)--(u1)--(u2)--(u3)--(u4)--(u5)--(u6)--(u7)--(u0);

\draw (w0) .. controls (wp1) and (wp2) .. (w3);
\draw (w1) .. controls (wp2) and (wp3) .. (w4);
\draw (w2) .. controls (wp3) and (wp4) .. (w5);
\draw (w3) .. controls (wp4) and (wp5) .. (w6);
\draw (w4) .. controls (wp5) and (wp6) .. (w7);
\draw (w5) .. controls (wp6) and (wp7) .. (w0);
\draw (w6) .. controls (wp7) and (wp0) .. (w1);
\draw (w7) .. controls (wp0) and (wp1) .. (w2);

\node [vertex,label=above:{\small $v_{04}$}] (v04) at (0,1) {};
\node [vertex,label=below right:{\small $v_{26}$}] (v26) at (1,0) {};
\node [vertex,label=below:{\small $v_{15}$}] (v15) at (0,-1) {};
\node [vertex,label=below:{\small $v_{37}$}] (v37) at (-1,0) {};
\draw (v0)--(v04)--(v4);
\draw (v2)--(v26)--(v6);
\draw (v3)--(v37)--(v7);
\draw (v1)--(v15)--(v5);

\node [vertex,label=above right:{\small $a$}] (a) at (0.25,0.25) {};
\node [vertex,label=below left:{\small $b$}] (b) at (-0.125,-0.125){};

\draw (v04)--(a)--(v26);
\draw (v15)--(b)--(v37);

\draw [line width=5pt, opacity=0.5, red, rounded corners] (u0.center)--(u1.center)--(u2.center)--(v2.center)--(v26.center)--(a.center)--(v04.center)--(v0.center)--cycle;

\node [rvertex] at (v04) {};
\node [rvertex] at (a) {};
\node [rvertex] at (v26) {};
\node [rvertex] at (v0) {};
\node [rvertex] at (u0) {};
\node [rvertex] at (u1) {};
\node [rvertex] at (u2) {};
\node [rvertex] at (v2) {};

%\draw (w0) -- (w3) -- (w6)--(w1)--(w4)--(w7)--(w2)--(w5)--(w0);

\end{tikzpicture}
\caption{Tutte's $8$-cage minus one edge}
\label{fig:partialtutte2}
\end{center}
\end{figure}

Now we consider yet another $8$-cycle, this time the cycle $E =(u_0,u_1,u_2,v_2,v_{26},a,v_{04},v_0)$ shown in red on \cref{fig:partialtutte2}. On this cycle, $a$ and $u_1$ are antipodal vertices, and the third neighbour of $u_1$ is $v_1$. So the third neighbour of $a$ off $E$ must share a common neighbour with $v_1$, i.e., it must be adjacent to a vertex that is distance $2$ from $v_1$. The only vertex with degree less than $3$ that meets that criterion is $b$, and so $a \sim b$ and $G$ is Tutte's $8$-cage.

\input girth9plus

%% file: girth3.tex
\subsection{Girth three}

In this subsection we show that the only cubic graph $G$ of girth $3$ with no eigenvalues in $(-2,0)$ is the $3$-prism $K_3 \mathbin{\square} K_2$.

% This is accomplished through a number of short lemmas where the main tactic is showing that certain subgraphs lead to negative principal minors of the matrix $M = A(A+2I)$, and therefore cannot occur in a graph with no eigenvalues in $(-2,0)$.

% \krystalsays{In fact, in the following, $G$ must have eigenvalue $-1$. The $M$ we have has eigenvalue $-1$, but we do not even need it because the end point of the edge on 2 triangles are closed twins. }

% \gordonsays{You are right, I have put in the simpler argument}

\begin{lemma}\label{lem:twotri}
If $G$ is a cubic graph of girth $3$ with two triangles sharing an edge, then $G$ has $-1$ as an eigenvalue.
\end{lemma}

\begin{proof}
Let $\{u,v,w\}$ and $\{u,v,x\}$ be two triangles in $G$ sharing the edge $\{u,v\}$. If $x$ is a vector with entries indexed by $V(G)$ whose only non-zero entries are $x_u=1$ and $x_v=-1$, then $Ax = -x$.
\end{proof}

% \begin{proof}
% Let $\{u,v,w\}$ and $\{u,v,x\}$ be two triangles in $G$. Then $u$ and $v$ are adjacent and have two common neighbours, and so $M_{uv} = 4$. Therefore 
% \[
% M_{\{u,v\},\{u,v\}} = \left[
% \begin{array}{cc}3&4\\4&3\end{array}
% \right]
% \]
% which has determinant $-7$.
% \end{proof}

%Here is Sym(\eps)ilon

\begin{lemma}\label{lem:g3-prism}
If $G$ is a cubic graph of girth $3$ with no 
eigenvalues in $(-2,0)$, then $G$ is the 
$3$-prism $K_3 \mathbin{\square} K_2$.

% If $G$ is a cubic graph of girth $3$ with no two triangles sharing an edge, then $G$ is either the $3$-prism $K_3 \mathbin{\square} K_2$ or $M(G)$ contains 
% \[
% M' =
% \left[
% \begin{array}{ccc}
% 3&3&1\\
% 3&3&2\\
% 1&2&3
% \end{array}
% \right]
% \]
% as a principal submatrix, and an eigenvalue in the interval $(-2,0)$.
\end{lemma}

\begin{proof}
Let $T = \{u,v,w\}$ be a triangle in $G$ and let $u'$, $v'$, $w'$ denote the third neighbours of $u$, $v$, $w$ respectively. By \cref{lem:twotri}, these vertices are distinct and so $G$ contains the configuration shown in the first figure of \cref{fig:tris}.

Suppose first that the three vertices 
$\{u',v',w'\}$ do not form a triangle, in which case we may assume without loss of generality that $v'$ and $w'$ are not adjacent. Then $v$ and $w'$ have a unique common neighbour and taking $S = \{v,w,w'\}$, we have 
\begin{equation}\label{mymat}
M_{SS} =
\kbordermatrix{
 &v&w&w'\\
v&3&3&1\\
w&3&3&2\\
w'&1&2&3
},
\end{equation}
which has determinant $-3$. Therefore $M$ is not positive semidefinite contradicting the hypotheses on $G$.

Otherwise, $\{u',v',w'\}$ do form a triangle, and the entire graph is the $3$-prism. 
\end{proof}

\begin{figure}
\begin{center}
\begin{tikzpicture}
\tikzstyle{vertex}=[fill=blue!25, circle, draw=black, inner sep=0.8mm]
\tikzstyle{rvertex}=[fill=red, circle, draw=black, inner sep=0.8mm]
\node [vertex,label={[label distance=-1mm]0:$u$}] (u) at (90:.8cm) {};
\node [vertex,label={[label distance=-1.5mm]135:$v$}] (v) at (210:0.8cm) {};
\node [vertex,label={[label distance=-1.5mm]30:$w$}] (w) at (330:0.8cm) {};
\node [vertex, label={[label distance=-1mm]0:$u'$}] (up) at (90:2cm) {};
\node  [vertex,label={[label distance=-1.5mm]135:$v'$}] (vp) at (210:2cm) {};
\node [vertex,label={[label distance=-1.5mm]30:$w'$}] (wp) at (330:2cm) {};

\node [vertex] (wp) at (330:2cm) {};
\draw (u)--(v)--(w)--(u);
%\draw (up)--(vp)--(wp)--(up);
\draw (u)--(up);
\draw (v)--(vp);
\draw (w)--(wp);

\pgftransformxshift{7cm};

\tikzstyle{vertex}=[fill=blue!25, circle, draw=black, inner sep=0.8mm]
\node [vertex,label={[label distance=-1mm]0:$u$}] (u) at (90:.8cm) {};
\node [vertex,label={[label distance=-1.5mm]135:$v$}] (v) at (210:0.8cm) {};
\node [vertex,label={[label distance=-1.5mm]30:$w$}] (w) at (330:0.8cm) {};
\node [vertex, label={[label distance=-1mm]0:$u'$}] (up) at (90:2cm) {};
\node  [vertex,label={[label distance=-1.5mm]135:$v'$}] (vp) at (210:2cm) {};
\node [vertex,label={[label distance=-1.5mm]30:$w'$}] (wp) at (330:2cm) {};

\node [vertex] (wp) at (330:2cm) {};
\draw (u)--(v)--(w)--(u);
%\draw (up)--(vp)--(wp)--(up);
\draw (u)--(up);
\draw (v)--(vp);
\draw (w)--(wp);

\node [rvertex] at (v) {};
\node [rvertex] at (w) {};
\node [rvertex] at (wp) {};

\draw [thick,dashed] (vp)--(wp);

\end{tikzpicture}
\end{center}
\caption{Configurations around a triangle}
\label{fig:tris}
\end{figure}

The ``quantitative version'' of this lemma is the following:

\begin{lemma}\label{lem:quant} Let $\tau'\approx -0.201912$ be the smallest root of $x^3 - 9x^2 + 13x + 3$. 
    If $G$ is a cubic graph of girth $3$, then either $G$ is isomorphic to the $3$-prism $K_3 \mathbin{\square} K_2$ or it has an eigenvalue in $[	-1 - \sqrt{1+ \tau'}, 	-1 + \sqrt{1+ \tau'}  ]
\approx [-1.893358, -0.106642]$. In particular, the only cubic graph of girth $3$ with no eigenvalues in $(-2,0)$ is the $3$-prism.
\end{lemma}

\begin{proof}
    If $G$ has two triangles sharing a edge, then $G$ has $-1$ as an eigenvalue, which is in the interval. Otherwise we follow the proof of \cref{lem:g3-prism} and obtain that $G$ is either isomorphic to the $3$-prism or $M(G)$ has the matrix shown in \cref{mymat} as a principal submatrix. The characteristic polynomial of this matrix is 
    $x^3 - 9x^2 + 13x + 3$
    which has smallest eigenvalue $\tau'\approx -0.201912$. The statement now follows from \cref{lem:mtau}. 
\end{proof}

%% file: girth4NEW.tex
\section{Girth four}

In this subsection we will prove that if $G$ is a cubic graph of girth $4$ with no eigenvalues in $(-2,0)$ then $G$ is either $K_{3,3}$ or $G$ is isomorphic to $X(n)$ for some $n \geqslant 2$.

We will need a sequence of lemmas progressively restricting the possible structure of $G$. So let $C = (u_0,u_1,u_2,u_3)$ be a $4$-cycle and let $v_i$ denote the unique neighbour of $u_i$ off $C$. The vertices $V = \{v_0, v_1, v_2, v_3\}$ may not be distinct but because $G$ has girth $4$ it follows that $v_i \not= v_{i \pm 1}$. So it may be the case that $v_0=v_2$ and $v_1 = v_3$, or (without loss of generality) $v_0 = v_2$ but $v_1 \not= v_3$, or that $v_0 \not= v_2$ and $v_1 \not= v_3$. These three possibilities are depicted in \cref{fig:corona4}: we note that the third graph is $K_{3,3}$ minus an edge.

\begin{figure}
\begin{center}
\begin{tikzpicture}
\tikzstyle{vertex}=[fill=blue!25, circle, draw=black, inner sep=0.7mm]

\node [vertex] (u0) at (3,3) {};
\node [vertex] (u1) at (0,3) {};
\node [vertex] (u2) at (0,0) {};
\node [vertex] (u3) at (3,0) {};

\node [vertex] (v0) at (2,2) {};
\node [vertex] (v1) at (1,2) {};
\node [vertex] (v2) at (1,1) {};
\node [vertex] (v3) at (2,1) {};

\node [above right] at (u0) {\small $u_0$};
\node [above left] at (u1) {\small $u_1$};
\node [below left] at (u2) {\small $u_2$};
\node [below right] at (u3) {\small $u_3$};

\node [right] at (v0) {\small $v_0$};
\node [left] at (v1) {\small $v_1$};
\node [left] at (v2) {\small $v_2$};
\node [right] at (v3) {\small $v_3$};

\draw (u0)--(u1)--(u2)--(u3)--(u0);
\draw (u0)--(v0);
\draw (u1)--(v1);
\draw (u2)--(v2);
\draw (u3)--(v3);

\node at (1.5,-0.75) {$C_1$};

\pgftransformxshift{5cm}
\node [vertex] (u0) at (3,3) {};
\node [vertex] (u1) at (0,3) {};
\node [vertex] (u2) at (0,0) {};
\node [vertex] (u3) at (3,0) {};

\node [vertex] (v0) at (2,2) {};
\node [vertex] (v1) at (1,2) {};
%\node [vertex] (v2) at (1,1) {};
\node [vertex] (v3) at (2,1) {};

\draw (u0)--(u1)--(u2)--(u3)--(u0);
\draw (u0)--(v0);
\draw (u1)--(v1);
\draw (u2)--(v0);
\draw (u3)--(v3);
\node [above right] at (u0) {\small $u_0$};
\node [above left] at (u1) {\small $u_1$};
\node [below left] at (u2) {\small $u_2$};
\node [below right] at (u3) {\small $u_3$};

\node [right] at (v0) {\small $v_0$};
\node [left] at (v1) {\small $v_1$};
%\node [left] at (v2) {\small $v_2$};
\node [right] at (v3) {\small $v_3$};

\node at (1.5,-0.75) {$C_2$};
\pgftransformxshift{5cm}
\node [vertex] (u0) at (3,3) {};
\node [vertex] (u1) at (0,3) {};
\node [vertex] (u2) at (0,0) {};
\node [vertex] (u3) at (3,0) {};

\node [vertex] (v0) at (2,2) {};
\node [vertex] (v1) at (1,2) {};
% \node [vertex] (v2) at (1,1) {};
% \node [vertex] (v3) at (2,1) {};

\draw (u0)--(u1)--(u2)--(u3)--(u0);
\draw (u0)--(v0);
\draw (u1)--(v1);
\draw (u2)--(v0);
\draw (u3)--(v1);
\node [above right] at (u0) {\small $u_0$};
\node [above left] at (u1) {\small $u_1$};
\node [below left] at (u2) {\small $u_2$};
\node [below right] at (u3) {\small $u_3$};

\node [right] at (v0) {\small $v_0$};
\node [left] at (v1) {\small $v_1$};
%\node [left] at (v2) {\small $v_2$};
%\node [right] at (v3) {\small $v_3$};
\node at (1.5,-0.75) {$C_3$};
\end{tikzpicture}

\end{center}
\caption{Corona of a $4$-cycle}
\label{fig:corona4}
\end{figure}

\begin{lemma}\label{lem:girth4-k33-edge}
Let $G$ be a cubic graph of girth $4$ with no eigenvalues in $(-2,0)$. If $G$ has a $4$-cycle whose corona is isomorphic to $C_3$ in \cref{fig:corona4}), then $G \cong K_{3,3}$.
\end{lemma}

\begin{proof}
Consider the distance $d_G(v_0,v_1)$ between $v_0$ and $v_1$ in $G$. If $d_G(v_0,v_1) = 1$, then  $G \cong K_{3,3}$. If $d_G(v_0,v_1) = 2$, then $v_0$ and $v_1$ have a common neighbour, say $v_{01}$, and $M_{v_0v_1} = 1$ and otherwise $d_G(v_0,v_1) > 2$ and so $M_{v_0v_1} = 0$. 

In the former case, let $T = \{u_0,u_3,v_0,v_1,v_{01}\}$ in which case
\[
M_{TT} = \kbordermatrix{
& u_0&u_3&v_0&v_1&v_{01}\\
u_0 & 3 & 2 & 2 & 2 & 1\\
u_3 & 2 & 3 & 2 & 2 & 1\\
v_0 & 2 & 2 & 3 & 1 & 2\\
v_1 & 2 & 2 & 1 & 3 & 2\\
v_{01} & 1 & 1 & 2 & 2 & 3
}
\]
which has determinant $-8$. In the latter case, let $T = \{u_0,u_1,v_0,v_1\}$ in which case
\[
M_{TT} = \kbordermatrix{
&u_0&u_1&v_0&v_1\\
u_0&3&2&2&2\\
u_1&2&3&2&2\\
v_0&2&2&3&0\\
v_1&2&2&0&3\\
}
\]
which has determinant $-3$. As $M(G)$ is positive semidefinite, neither of these cases can occur.
\end{proof}

% \begin{lemma}\label{lem:girth4-k33-edge}
% Let $G$ be a cubic 
% If $G$ is a cubic graph of girth 4 with no eigenvlaues in $(0,2)$that contains a subgraph isomorphic to $K_{3,3}-e$ (where $e$ is an arbitrary edge of $K_{3,3}$) then $G = K_{3,3}$ or $M(G)$ contains either $M_1$ or $M_2$ as a principal submatrix. 
% \end{lemma}

% \begin{lemma}\label{lem:girth4-k33-edge}
% If $G$ is a cubic graph of girth 4 that contains a subgraph isomorphic to $K_{3,3}-e$ (where $e$ is an arbitrary edge of $K_{3,3}$) then $G = K_{3,3}$ or $M(G)$ contains either $M_1$ or $M_2$ as a principal submatrix. 
% \end{lemma}

% \begin{proof}
% Suppose that $G$ contains $K_{3,3}-e$ as depicted in the first diagram of \cref{fig:k33e}. The distance from $u$ to $x$ is either $1$, $2$ or more than $2$. In the first case, $u \sim x$ and so $G = K_{3,3}$ which is a genuine example of a graph with no eigenvalues in $(-2,0)$. In the second case, $u$ and $x$ have a common neighbour $y$, and setting 
% $S = \{u,v,w,x,y\}$ we get
% \[
% M_{SS} = 
% \kbordermatrix{
% & u&v&w&x&y\\
% u&3&2&2&1&2\\
% v&2&3&2&2&1\\
% w&2&2&3&2&1\\
% x&1&2&2&3&2\\
% y&2&1&1&2&3\\
% }
% \]
% which has determinant $-8$ and is equal to $M_1$. 

% In the third case $M_{ux} = 0$ and if we take $T =  \{u,v,w,x\}$ then
% \[
% M_{TT} = 
% \kbordermatrix{
% &u&v&w&x\\
% u&3&2&2&0\\
% v&2&3&2&2\\
% w&2&2&3&2\\
% x&0&2&2&3
% }
% \]
% which is equal to $M_2$.
% \end{proof}

\begin{figure}
\begin{center}
\begin{tikzpicture}[scale=1.25]
\tikzstyle{vertex}=[fill=blue!25, circle, draw=black, inner sep=0.8mm]
\node [vertex] (u0) at (0,1.5) {};
\node [vertex] (u2) at (0,0.5) {};
\node [vertex] (u1) at (1.5,2) {};
\node [vertex] (u3) at (1.5,1) {};
\node [vertex] (v0) at (1.5,0) {};
\node [vertex] (v1) at (2.5,2) {};
\node [vertex] (v3) at (2.5,1) {};
\node [vertex] (v0p) at (2.5,0) {};
\node [left] at (u0) {$u_0$\,};
\node [left] at (u2) {$u_2$\,};
\node [above right] at (u1) {$u_1$};
\node [above right] at (u3) {$u_3$};
\node [above right] at (v0) {$v_0$};
\node [above right] at (v1) {$v_1$};
\node [above right] at (v3) {$v_3$};
\node [above right] at (v0p) {$v_0'$};
\draw (u0)--(u1);
\draw (u0)--(u3);
\draw (u0)--(v0);
\draw (u2)--(u1)--(v1);
\draw (u2)--(u3)--(v3);
\draw (u2)--(v0)--(v0p);

% \node [right] (t) at (3.5,1.5) {\small not adjacent};
% \draw [->](t)--(v1);
% \draw [->] (t)--(v3);
\end{tikzpicture}
\hspace{2cm}
\begin{tikzpicture}[scale=1.25]
\tikzstyle{vertex}=[fill=blue!25, circle, draw=black, inner sep=0.8mm]
\tikzstyle{rvertex}=[fill=red, circle, draw=black, inner sep=0.8mm]
\node [vertex] (u0) at (0,1.5) {};
\node [vertex] (u2) at (0,0.5) {};
\node [vertex] (u1) at (1.5,2) {};
\node [vertex] (u3) at (1.5,1) {};
\node [vertex] (v0) at (1.5,0) {};
\node [vertex] (v1) at (2.5,2) {};
\node [vertex] (v3) at (2.5,1) {};
\node [vertex] (v0p) at (2.5,0) {};
\node [left] at (u0) {$u_0$\,};
\node [left] at (u2) {$u_2$\,};
\node [above right] at (u1) {$u_1$};
\node [above right] at (u3) {$u_3$};
\node [above right] at (v0) {$v_0$};
\node [above right] at (v1) {$v_1$};
\node [above right] at (v3) {$v_3$};
\node [above right] at (v0p) {$v_0'$};
\draw (u0)--(u1);
\draw (u0)--(u3);
\draw (u0)--(v0);
\draw (u2)--(u1)--(v1);
\draw (u2)--(u3)--(v3);
\draw (u2)--(v0)--(v0p);

\draw [thick, dotted] (v1)--(v3);
\node [right] (tx) at (3.2,1.5) {\small non-edge};
\draw [->] (tx)--(2.6,1.5);

\node [rvertex] at (u1) {};
\node [rvertex] at (v1) {};
\node [rvertex] at (u3) {};
\node [rvertex] at (v3) {};

\draw [ultra thick, red] (u1)--(v1);
\draw [ultra thick, red] (u3)--(v3);

% \node [right] (t) at (3.5,1.5) {\small not adjacent};
% \draw [->](t)--(v1);
% \draw [->] (t)--(v3);
\end{tikzpicture}

\caption{Configuration around a corona isomorphic to $C_2$.}
\label{fig:corc2}
\end{center}
\end{figure}

\begin{lemma}
Let $G$ be a cubic graph of girth $4$ with no eigenvalues in $(-2,0)$. Then $G$ has no $4$-cycle whose corona is isomorphic to $C_2$ in \cref{fig:corona4}.
\end{lemma}

\begin{proof}
Suppose for a contradiction that $G$ does have such a $4$-cycle and that the corresponding corona is labelled as in \cref{fig:corona4}. Let $v_0'$ be the third neighbour of $v_0$. If $v_0' = v_1$, then $\{u_0,u_1,u_2,v_0\}$ is a $4$-cycle with a corona isomorphic to $C_3$ which is not possible. So $\{v_1,v_3,v_0'\}$ are distinct vertices and we have the configuration depicted in \cref{fig:corc2} (redrawn to highlight the symmetry of the situation). As $G$ has girth $4$ it cannot be the case that the subgraph induced by $\{v_1,v_3,v_0'\}$ is a triangle and so we may assume without loss of generality that $v_1 \not\sim v_3$. Thus we have the configuration depicted in \cref{fig:corc2}.
Taking $T = \{u_1,v_1,u_3,v_3\}$ we get
\[
M_{TT} = 
\kbordermatrix{
&u_1&v_1&u_3&v_3\\
u_1&3&2&2&0\\
v_1&2&3&0&\alpha\\
u_3&2&0&3&2\\
v_3&0&\alpha&2&3
}
\]
where $\alpha = 1$ if $d_G(v_1,v_3) = 2$ and $\alpha = 0$ if $d_G(v_1,v_3) > 2$. Now $\det M_{TT} = -11-16\alpha-5\alpha^2$ which is negative for $\alpha \in \{0,1\}$, thus contradicting the assumption that $M(G)$ is positive semidefinite.
\end{proof}

After these two lemmas, we know that the corona of every $4$-cycle is isomorphic to $C_1 \cong C_4 \circ K_1$. Next we consider the possible edges and $2$-paths between the vertices $\{v_0,v_1,v_2,v_3\}$.

\begin{lemma}\label{c1lemma}
If $G$ is a cubic graph of girth $4$ with no eigenvalues in $(-2,0)$ and $G$ is not isomorphic to $K_{3,3}$ then in the corona of every $4$-cycle, labelled as $C_1$ in \cref{fig:corona4}),it must be the case that $v_0 \sim v_2$ and $v_1 \sim v_3$ and there are no other edges between pairs of vertices in $V$.
\end{lemma}

\begin{proof}
Suppose for a contradiction that $v_0 \not\sim v_2$. Taking $T = \{u_0,v_0,u_2,v_2\}$ we get
\[
M_{TT} = 
\kbordermatrix{
&u_0&v_0&u_2&v_2\\
u_0&3&2&2&0\\
v_0&2&3&0&\alpha\\
u_2&2&0&3&2\\
v_2&0&\alpha&2&3
}
\]
where $\alpha \in \{0,1,2\}$, depending on how many common neighbours are shared by $v_0$ and $v_2$.  
As previously, $\det M_{TT} = -11-16\alpha-5\alpha^2$ which is negative for $\alpha \in \{0,1,2\}$, thus contradicting the assumption that $M(G)$ is positive semidefinite.

So suppose that $v_0 \sim v_2$ and $v_1 \sim v_3$ leading to the first configuration in \cref{fig:corc1}. Suppose for a contradiction that there is another edge, say $v_0v_1$ connecting two vertices in $V$, as depicted in the second diagram in \cref{fig:corc1}.

If we now take $T = \{u_0,u_1,u_2,v_2,v_0\}$ we have
\[
M_{TT} = 
\kbordermatrix{
&u_0&u_1&u_2&v_2&v_0\\
u_0&3&2&2&1&2\\
u_1&2&3&2&1&2\\
u_2&2&2&3&2&1\\
v_2&1&1&2&3&2\\
v_0&2&2&1&2&3
}
\]
and $\det M_{TT} = -8$, thereby contradicting the hypothesis that $M(G)$ is positive semidefinite.

Therefore the subgraph induced by the vertices of the corona of any $4$-cycle is isomorphic to $C_1$ from \cref{fig:corona4}.
\end{proof}

We reiterate that determining the entries of $M_{TT}$ is not entirely trivial, because it is always necessary consider the possibility that there may be extra edges between vertices in $T$ and/or $2$-paths between vertices in $T$ whose middle vertex is outside $T$.

\begin{figure}
\begin{center}
\begin{tikzpicture}
\tikzstyle{vertex}=[fill=blue!25, circle, draw=black, inner sep=0.7mm]
\tikzstyle{rvertex}=[fill=red, circle, draw=black, inner sep=0.8mm]

\node [vertex] (u0) at (3,3) {};
\node [vertex] (u1) at (0,3) {};
\node [vertex] (u2) at (0,0) {};
\node [vertex] (u3) at (3,0) {};

\node [vertex] (v0) at (2,2) {};
\node [vertex] (v1) at (1,2) {};
\node [vertex] (v2) at (1,1) {};
\node [vertex] (v3) at (2,1) {};

\node [above right] at (u0) {\small $u_0$};
\node [above left] at (u1) {\small $u_1$};
\node [below left] at (u2) {\small $u_2$};
\node [below right] at (u3) {\small $u_3$};

\node [right] at (v0) {\small $v_0$};
\node [left] at (v1) {\small $v_1$};
\node [left] at (v2) {\small $v_2$};
\node [right] at (v3) {\small $v_3$};

\draw (u0)--(u1)--(u2)--(u3)--(u0);
\draw (u0)--(v0);
\draw (u1)--(v1);
\draw (u2)--(v2);
\draw (u3)--(v3);

\draw (v0)--(v2);
\draw (v1)--(v3);

\pgftransformxshift{6cm}

\node [vertex] (u0) at (3,3) {};
\node [vertex] (u1) at (0,3) {};
\node [vertex] (u2) at (0,0) {};
\node [vertex] (u3) at (3,0) {};

\node [vertex] (v0) at (2,2) {};
\node [vertex] (v1) at (1,2) {};
\node [vertex] (v2) at (1,1) {};
\node [vertex] (v3) at (2,1) {};

\node [above right] at (u0) {\small $u_0$};
\node [above left] at (u1) {\small $u_1$};
\node [below left] at (u2) {\small $u_2$};
\node [below right] at (u3) {\small $u_3$};

\node [right] at (v0) {\small $v_0$};
\node [left] at (v1) {\small $v_1$};
\node [left] at (v2) {\small $v_2$};
\node [right] at (v3) {\small $v_3$};

\draw (u0)--(u1)--(u2)--(u3)--(u0);
\draw (u0)--(v0);
\draw (u1)--(v1);
\draw (u2)--(v2);
\draw (u3)--(v3);

\draw (v0)--(v2);
\draw (v1)--(v3);
\draw (v0)--(v1);

\node [rvertex] at (u0) {};
\node [rvertex] at (u1) {};
\node [rvertex] at (u2) {};
\node [rvertex] at (v2) {};
\node [rvertex] at (v0) {};

\draw [ultra thick, red] (u0)--(u1)--(u2)--(v2)--(v0)--(u0);

\end{tikzpicture}
\end{center}

\caption{Graph induced by corona $C_1$.}
\label{fig:corc1}
\end{figure}

Now we consider how a copy of $C_1$ can be extended in a cubic graph of girth $4$ with no eigenvalues in $(-2,0)$, and show that the only possibility is a cyclic sequence of alternating $4$-cycles and single edges connected as depicted in \cref{fig:xk5}.

\begin{theorem}
If $G$ is a cubic graph of girth $4$ with no eigenvalues in $(-2,0)$, then $G \cong K_{3,3}$ or $G \cong X(n)$ for some $n \geqslant 2$.
\end{theorem}

\begin{proof}
If $G \cong K_{3,3}$ we are done. Otherwise, $G$ has an induced subgraph isomorphic to the $6$-vertex gadget of \cref{fig:gadget}. As $(v_2,v_4,v_3,v_5)$ is a $4$-cycle, its corona must be isomorphic to $C_1$, and so there are adjacent vertices $w_0$, $w_1$ such that $w_0 \sim v_4$ and $w_1 \sim v_5$. 

By \cref{c1lemma}, these new vertices $w_0$, $w_1$ cannot be adjacent to either $v_0$ or $v_1$. So $w_0$ is adjacent to a new vertex $w_2$, and $w_1$ is adjacent to a new vertex $w_3$. The vertices $w_2$ and $w_3$ are distinct and must be non-adjacent (to avoid creating a $4$-cycle with corona not isomorphic to $C_1$), leading to the graph shown in the first diagram of \cref{fig:xtendgadget}.  Next we wish to show that there are no edges between $\{v_0,v_1\}$ and $\{w_2,w_3\}$.  If there were such an edge, say $v_0 \sim w_2$, then taking $T = \{v_0,v_1,v_2,v_3,w_0\}$ we would have
\[
M_{TT} =
\kbordermatrix{
&v_0&v_1&v_2&v_3&w_0\\
v_0& 3&2&2&1&1 \\
v_1& 2&3&1&2&0\\
v_2&2&1&3&2&1\\
v_3& 1&2&2&3&1\\
w_0& 1&0&1&1&3
}
\]
for which $\det M_{TT} = -8$.  Therefore $w_2$ is adjacent to two new vertices, which we denote $w_4$ and $w_5$, yielding the graph in the second diagram of \cref{fig:xtendgadget}. All that remains is to show that $w_3$ is also adjacent to these two new vertices, thereby demonstrating that the gadget repeats. If $w_3$ was not adjacent to one of $\{w_4,w_5\}$, say $w_4$, then taking 
$T = \{v_4,v_5,w_0,w_1,w_4\}$ we would obtain
\[
M_{TT} = 
\kbordermatrix{
&v_4&v_5,&w_0&w_1&w_4\\
v_4&3&2&2&1&0\\
v_5&2&3&1&2&0\\
w_0&2&1&3&2&1\\
w_1&1&2&2&3&0\\
w_4&0&0&1&0&3
}
\]
for which $\det M_{TT} = -8$. Therefore $w_3 \sim w_4$ and by symmetry $w_3 \sim w_5$ and the entire $6$-vertex gadget is replicated, as in the third graph of \cref{fig:xtendgadget}.  Now $(w_2,w_4,w_3,w_5)$ is a new $4$-cycle and so the third neighbours of $w_4$ and $w_5$ must be adjacent. These neighbours can either be $v_0$ and $v_1$, in which case the graph is cubic and isomorphic to $X(2)$, or they can be two new vertices, in which case the argument can be repeated to add another gadget to the growing graph. After each six-vertex gadget is added, the edges from the two end-vertices may either ``wrap round'' to $v_0$ and $v_1$, creating a cubic graph isomorphic to $X(n)$ for some $n$, or they may connect to two new vertices, which must be adjacent and form the first two vertices of the next six-vertex gadget.
\end{proof}

\begin{figure}
\begin{center}
\begin{tikzpicture}
\tikzstyle{vertex}=[fill=blue!25, circle, draw=black, inner sep=0.8mm]
\node [vertex] (v0) at (0,0) {};
\node [vertex] (v1) at (0,-1) {};
\node [vertex] (v2) at (1,0) {};
\node [vertex] (v3) at (1,-1) {};
\node [vertex] (v4) at (2,0) {};
\node [vertex] (v5) at (2,-1) {};

\node [vertex] (w0) at (3,0) {};
\node [vertex] (w1) at (3,-1) {};

\node [vertex] (w2) at (4,0) {};
\node [vertex] (w3) at (4,-1) {};

\node [above] at (v0.north) {\small $v_0$};
\node [below] at (v1.south) {\small $v_1$};
\node [above] at (v2.north) {\small $v_2$};
\node [below] at (v3.south) {\small $v_3$};
\node [above] at (v4.north) {\small $v_4$};
\node [below] at (v5.south) {\small $v_5$};
\node [above] at (w0.north) {\small $w_0$};
\node [below] at (w1.south) {\small $w_1$};
\node [above] at (w2.north) {\small $w_2$};
\node [below] at (w3.south) {\small $w_3$};
\draw (v0)--(v2);
\draw (v4)--(v3);
\draw (v5)--(v2);
\draw (v1)--(v3);
\draw (v2)--(v4);
\draw (v3)--(v5);
\draw (v0)--(v1);
\draw (v1)--(-0.65,-1);
\draw (v0)--(-0.65,0);

\draw (v4)--(w0)--(w1)--(v5);
\draw (w0)--(w2);
\draw (w1)--(w3);

% \draw (v4)--(2.65,0);
% \draw (v5)--(2.65,-1);

\end{tikzpicture}
\hspace{1cm}
\begin{tikzpicture}
\tikzstyle{vertex}=[fill=blue!25, circle, draw=black, inner sep=0.8mm]
\node [vertex] (v0) at (0,0) {};
\node [vertex] (v1) at (0,-1) {};
\node [vertex] (v2) at (1,0) {};
\node [vertex] (v3) at (1,-1) {};
\node [vertex] (v4) at (2,0) {};
\node [vertex] (v5) at (2,-1) {};

\node [vertex] (w0) at (3,0) {};
\node [vertex] (w1) at (3,-1) {};

\node [vertex] (w2) at (4,0) {};
\node [vertex] (w3) at (4,-1) {};

\node [vertex] (w4) at (5,0) {};
\node [vertex] (w5) at (5,-1) {};

\node [above] at (v0.north) {\small $v_0$};
\node [below] at (v1.south) {\small $v_1$};
\node [above] at (v2.north) {\small $v_2$};
\node [below] at (v3.south) {\small $v_3$};
\node [above] at (v4.north) {\small $v_4$};
\node [below] at (v5.south) {\small $v_5$};
\node [above] at (w0.north) {\small $w_0$};
\node [below] at (w1.south) {\small $w_1$};
\node [above] at (w2.north) {\small $w_2$};
\node [below] at (w3.south) {\small $w_3$};

\node [above] at (w4.north) {\small $w_4$};
\node [below] at (w5.south) {\small $w_5$};

\draw (v0)--(v2);
\draw (v4)--(v3);
\draw (v5)--(v2);
\draw (v1)--(v3);
\draw (v2)--(v4);
\draw (v3)--(v5);
\draw (v0)--(v1);
\draw (v1)--(-0.65,-1);
\draw (v0)--(-0.65,0);

\draw (v4)--(w0)--(w1)--(v5);
\draw (w0)--(w2);
\draw (w1)--(w3);

\draw (w2)--(w4);
\draw (w2)--(w5);

% \draw (v4)--(2.65,0);
% \draw (v5)--(2.65,-1);

\end{tikzpicture}

\vspace{0.5cm}

\begin{tikzpicture}
\tikzstyle{vertex}=[fill=blue!25, circle, draw=black, inner sep=0.8mm]
\node [vertex] (v0) at (0,0) {};
\node [vertex] (v1) at (0,-1) {};
\node [vertex] (v2) at (1,0) {};
\node [vertex] (v3) at (1,-1) {};
\node [vertex] (v4) at (2,0) {};
\node [vertex] (v5) at (2,-1) {};

\node [vertex] (w0) at (3,0) {};
\node [vertex] (w1) at (3,-1) {};

\node [vertex] (w2) at (4,0) {};
\node [vertex] (w3) at (4,-1) {};

\node [vertex] (w4) at (5,0) {};
\node [vertex] (w5) at (5,-1) {};

\node [above] at (v0.north) {\small $v_0$};
\node [below] at (v1.south) {\small $v_1$};
\node [above] at (v2.north) {\small $v_2$};
\node [below] at (v3.south) {\small $v_3$};
\node [above] at (v4.north) {\small $v_4$};
\node [below] at (v5.south) {\small $v_5$};
\node [above] at (w0.north) {\small $w_0$};
\node [below] at (w1.south) {\small $w_1$};
\node [above] at (w2.north) {\small $w_2$};
\node [below] at (w3.south) {\small $w_3$};

\node [above] at (w4.north) {\small $w_4$};
\node [below] at (w5.south) {\small $w_5$};

\draw (v0)--(v2);
\draw (v4)--(v3);
\draw (v5)--(v2);
\draw (v1)--(v3);
\draw (v2)--(v4);
\draw (v3)--(v5);
\draw (v0)--(v1);
\draw (v1)--(-0.65,-1);
\draw (v0)--(-0.65,0);

\draw (v4)--(w0)--(w1)--(v5);
\draw (w0)--(w2);
\draw (w1)--(w3);

\draw (w2)--(w4);
\draw (w2)--(w5);

\draw (w3)--(w4);
\draw (w3)--(w5);

% \draw (v4)--(2.65,0);
% \draw (v5)--(2.65,-1);

\end{tikzpicture}
\caption{Extending the gadget}
\label{fig:xtendgadget}
\end{center}
\end{figure}

%% file: girth5X.tex
\section{Girth five}

This section is dedicated to proving the following:

\begin{theorem}\label{thm:girth5}
If $G$ is a cubic graph of girth $5$ with no eigenvalues in $(-2,0)$ then $G$ is either the Petersen graph or the dodecahedron.
\end{theorem}

\newcommand\son[1]{\ensuremath N_{\geq 2}(#1)}

If $S$ is a set of vertices in a graph, then we define the \emph{strong open neighbourhood} of $S$ by 
\[
\son{S} = \{ v \in V(G) \backslash S \mathbin{:} |N(v) \cap S| \geqslant 2\}.
\]
In other words the strong open neighbourhood of $S$ consists of the vertices outside $S$ that are adjacent to at least two vertices inside $S$.

For $S \subseteq V(G)$, let $\widetilde{S}$ denote the subgraph of $G$ with
\begin{align*}
V(\widetilde{S}) &= S \cup \son{S},\\
E(\widetilde{S}) &= \{ \{u,v\} \in E(G) : u \in S, v \in S \cup \son{S}\}.
\end{align*}
The purpose of this definition is that $\widetilde{S}$ is the smallest subgraph of $G$ from which $M_{SS}$ can be determined. We note here that $\widetilde{S}$ is not necessarily an induced subgraph of $G$, because it does not contain any edges with both end-vertices in $\son{S}$ (these edges do not contribute to $M_{SS}$).

\begin{figure}
\begin{center}
\begin{tikzpicture}
\tikzstyle{vertex}=[fill=blue!25, circle, draw=black,inner sep=0.7mm]
\foreach \x in {0,1,...,4} {
  \node [vertex] (u\x) at (90+72*\x:2.5cm) {};
  \node at (90+72*\x:2.9cm) {$u_\x$};
    \node [vertex] (v\x) at (90+72*\x:1.5cm) {};
  \node at (90+72*\x:1.1cm) {$v_\x$};
  \draw (u\x)--(v\x);
}

% \node [vertex] (v0) at (90 {};
% \node [vertex] (v1) at (-1.902113,0.618034) {};
% \node [vertex] (v2) at (-1.175571,-1.618034) {};
% \node [vertex] (v3) at (1.175571,-1.618034) {};
% \node [vertex] (v4) at (1.902113,0.618034) {};
% \node [vertex] (v5) at (0.000000,1.000000) {};
% \node [vertex] (v6) at (-0.951057,0.309017) {};
% \node [vertex] (v7) at (-0.587785,-0.809017) {};
% \node [vertex] (v8) at (0.587785,-0.809017) {};
% \node [vertex] (v9) at (0.951057,0.309017) {};
\draw (u0)--(u1)--(u2)--(u3)--(u4)--(u0);
%\node at (0,-2.25) {$X_{0}$};
\end{tikzpicture}
\end{center}
\caption{The corona product $C_5 \circ K_1$}
\label{fig:corona5}
\end{figure}

If $C$ is a $5$-cycle in $G$, then its neighbours are distinct and $\cor{C}$ is a subgraph (not necessarily induced) isomorphic to $C_5 \circ K_1$ (see \cref{fig:corona5}). If $G$ has no eigenvalues in $(-2,0)$ then, setting $S = V(\cor{C})$, it must be the case that $M_{SS}$ is positive semi-definite. We shall see that this places severe constraints on the possibilities for the subgraph $\widetilde{S}$.

Starting with the vertices and edges of $C_5 \circ K_1$, we add new edges (between vertices of $C_5 \circ K_1$)) and new vertices (adjacent to $2$ or $3$ vertices of $C_5 \circ K_1$) in all possible ways that respect the degree and girth constraints. This leaves a long list of graphs as possibilities for $\widetilde{S}$, but the requirement that $M_{SS}$ be positive semidefinite means that many of these subgraphs simply cannot occur in a graph with no eigenvalues in $(-2,0)$. A short computation shows that  all but $13$ choices for $\widetilde{S}$ are ruled out by this requirement: these are depicted as $X_0$, $X_1$, $\ldots$, $X_{12}$ in the diagram \cref{fig:possibles}.

\begin{figure}
\begin{center}
\input X0
\hspace{0.95cm}
\input X1
\hspace{0.95cm}
\input X2

\vspace{0.275cm}

\input X3
\hspace{0.95cm}
\input X4
\hspace{0.95cm}
\input X5

\vspace{0.275cm}

\input X6
\hspace{0.95cm}
\input X7
\hspace{0.95cm}
\input X8

\vspace{0.275cm}

\input X9rot
\hspace{0.95cm}
\input X10
\hspace{0.95cm}
\input X11rot

\vspace{0.275cm}

\input X12
\end{center}
\caption{Possibilities for the subgraph $\widetilde{S}$ when $S = V(\cor{C_5})$.}
\label{fig:possibles}
\end{figure}

Looking ahead, the only graphs $\widetilde{S}$ that will actually occur in the list of graphs in \cref{mainresult} are $X_1$ when $G$ is the Petersen graph and $X_7$ when $G$ is the dodecahedron. In a few short lemmas, we rule out all of the other possibilities. Recall that $G$ may still contain edges or 2-paths that are not explicitly shown in $X_i$, but only if they connect two ``new vertices'' i.e., the vertices added after the initial $10$ vertices.

\begin{lemma}\label{notx0x2x3x4x8x10x12}
The subgraph $\widetilde{S}$ does not have a vertex of degree one and therefore $\widetilde{S} \not\cong X_0$, $X_2$, $X_3$, $X_4$, $X_8$, $X_{10}$, $X_{12}$.
\end{lemma}

\begin{proof}
Without loss of generality, suppose that $v_0$ has degree one in $\widetilde{C}$ and that $w_0$, $w_1 \ne u_0$ are the two neighbours of $v_0$ not in $\widetilde{C}$. Note that $w_0$ and $w_1$ are not adjacent to any of the vertices $\{v_1, v_2, v_3, v_4\}$ and therefore $M_{wu} = 0$ for all $w\in \{w_0,w_1\}$ and $v \in \{u_1, u_2, u_3, u_4\}$. Therefore if we take
\[
T = \{u_0, u_1, u_2, u_3, u_4, v_0, w_0, w_1\}
\]
it follows that 
\[
M_{TT} = 
\kbordermatrix{
      & u_0 & u_1 & u_2 & u_3 & u_4 & v_0 & w_0 & w_1 \\
u_0   & 3   & 2   & 1   & 1   & 2   & 2  & 1 & 1\\
u_1   & 2   & 3   & 2   & 1   & 1   & 1  & 0 & 0 \\
u_2   & 1   & 2   & 3   & 2   & 1   & 0  & 0 & 0\\
u_3   & 1   & 1   & 2   & 3   & 2 & 0 & 0 & 0 \\
u_4   & 2   & 1   & 1   & 2 & 3   & 1 & 0 & 0  \\
v_0   & 2   & 1   & 0   & 0 & 1   & 3  & 2 & 2  \\
w_0 & 1 & 0 & 0 & 0 & 0 & 2 & 3 & 1\\
w_1 & 1 & 0 & 0 & 0 & 0 & 2 & 1 & 3\\
}
\]
which has determinant $-4$.

\end{proof}

\begin{figure}
\begin{center}
 \begin{tikzpicture}[scale=1.1]
\tikzstyle{vertex}=[fill=blue!25, circle, draw=black,inner sep=0.7mm]
\tikzstyle{redvertex}=[fill=red, circle, draw=black,inner sep=0.7mm]
\node [vertex] (v0) at (0.000000,2.000000) {};
\node [vertex] (v1) at (-1.902113,0.618034) {};
\node [vertex] (v2) at (-1.175571,-1.618034) {};
\node [vertex] (v3) at (1.175571,-1.618034) {};
\node [vertex] (v4) at (1.902113,0.618034) {};
\node [vertex] (v5) at (0.000000,1.000000) {};
\node [vertex] (v6) at (-0.951057,0.309017) {};
\node [vertex] (v7) at (-0.587785,-0.809017) {};
\node [vertex] (v8) at (0.587785,-0.809017) {};
\node [vertex] (v9) at (0.951057,0.309017) {};
\node [vertex] (v10) at (-0.475528,0.654508) {};
\node [vertex] (v11) at (0.475528,0.654508) {};
\node [vertex] (v12) at (0.000000,-0.809017) {};
\node [above left] at (v10) {\small $v_{01}$};
\node [above right] at (v11) {\small $v_{04}$};
\draw (v0)--(v1);
\draw (v0)--(v4);
\draw (v0)--(v5);
\draw (v1)--(v2);
\draw (v1)--(v6);
\draw (v2)--(v3);
\draw (v2)--(v7);
\draw (v3)--(v4);
\draw (v3)--(v8);
\draw (v4)--(v9);
\draw (v5)--(v10);
\draw (v5)--(v11);
\draw (v6)--(v10);
\draw (v7)--(v12);
\draw (v8)--(v12);
\draw (v9)--(v11);
\node at (0,-2.25) {$X_{5}$};
\node [redvertex] at (v2) {};
\node [redvertex] at (v3) {};
\node [redvertex] at (v6) {};
\node [redvertex] at (v8) {};
\node [redvertex] at (v9) {};
\node [redvertex] at (v10) {};
\node [redvertex] at (v11) {};
\draw [ultra thick, red] (v6)--(v10);
\draw [ultra thick, red] (v11)--(v9);
\draw [ultra thick, red] (v2)--(v3)--(v8);
\end{tikzpicture}
\hspace{1cm}
\begin{tikzpicture}[scale=1.1]
\tikzstyle{vertex}=[fill=blue!25, circle, draw=black,inner sep=0.7mm]
\tikzstyle{redvertex}=[fill=red, circle, draw=black,inner sep=0.7mm]
\node [vertex] (v0) at (0.000000,2.000000) {};
\node [vertex] (v1) at (-1.902113,0.618034) {};
\node [vertex] (v2) at (-1.175571,-1.618034) {};
\node [vertex] (v3) at (1.175571,-1.618034) {};
\node [vertex] (v4) at (1.902113,0.618034) {};
\node [vertex] (v5) at (0.000000,1.000000) {};
\node [vertex] (v6) at (-0.951057,0.309017) {};
\node [vertex] (v7) at (-0.587785,-0.809017) {};
\node [vertex] (v8) at (0.587785,-0.809017) {};
\node [vertex] (v9) at (0.951057,0.309017) {};
\node [vertex] (v10) at (-0.475528,0.654508) {};
\node [vertex] (v11) at (0.475528,0.654508) {};
\node [vertex] (v12) at (-0.769421,-0.250000) {};
\node [vertex] (v13) at (0.000000,-0.809017) {};
\node [above left] at (v10) {\small $v_{01}$};
\node [above right] at (v11) {\small $v_{04}$};
\draw (v0)--(v1);
\draw (v0)--(v4);
\draw (v0)--(v5);
\draw (v1)--(v2);
\draw (v1)--(v6);
\draw (v2)--(v3);
\draw (v2)--(v7);
\draw (v3)--(v4);
\draw (v3)--(v8);
\draw (v4)--(v9);
\draw (v5)--(v10);
\draw (v5)--(v11);
\draw (v6)--(v10);
\draw (v6)--(v12);
\draw (v7)--(v12);
\draw (v7)--(v13);
\draw (v8)--(v13);
\draw (v9)--(v11);
\node [redvertex] at (v2) {};
\node [redvertex] at (v3) {};
\node [redvertex] at (v6) {};
\node [redvertex] at (v8) {};
\node [redvertex] at (v9) {};
\node [redvertex] at (v10) {};
\node [redvertex] at (v11) {};
\draw [ultra thick, red] (v6)--(v10);
\draw [ultra thick, red] (v11)--(v9);
\draw [ultra thick, red] (v2)--(v3)--(v8);
\node at (0,-2.25) {$X_{6}$};
\end{tikzpicture}
\end{center}
\caption{Ruling out $X_5$ and $X_6$}
\label{fig:x5x6}

\end{figure}

\begin{lemma}\label{notx5x6}
The subgraph $\widetilde{C} \not= X_5$, $X_6$.
\end{lemma}

\begin{proof}
In each case, let 
\[
T = \{u_2, u_3, v_1, v_3, v_{01}, v_4, v_{04}\}
\]
where $v_{01}$ is the common neighbour of $v_0$ and $v_1$, and $v_{04}$ is the common neighbour of $v_0$ and $v_4$, as depicted in \cref{fig:x5x6}. The only vertices in $T$ that are outside $\cor{C}$ are $v_{01}$ and $v_{04}$ but there can be no edge or 2-path between these two vertices by the girth constraint. Therefore 
\[
M_{TT} = 
\kbordermatrix{
& u_2& u_3& v_1& v_3& v_{01}& v_4& v_{04}\\
u_2 &3 &2 &1 &1 &0 &0 &0 \\
u_3 &2 &3 &0 &2 &0 &1 &0 \\
v_1 &1 &0 &3 &0 &2 &0 &0 \\
v_3 &1 &2 &0 &3 &0 &0 &0 \\
v_{01}&0 &0 &2 &0 &3 &0 &1 \\
v_{4}&0 &1 &0 &0 &0 &3 &2 \\
v_{04}&0 &0 &0 &0 &1 &2 &3 \\
}
\]
which has determinant $-36$ and therefore $M$ is not positive semidefinite.
\end{proof}

\begin{figure}
\begin{center}
\begin{tikzpicture}[scale=1.1]
\tikzstyle{vertex}=[fill=blue!25, circle, draw=black,inner sep=0.7mm]
\tikzstyle{redvertex}=[fill=red, circle, draw=black,inner sep=0.7mm]
\node [vertex] (v0) at (-1.902113,0.618034) {};
\node [vertex] (v1) at (-1.175571,-1.618034) {};
\node [vertex] (v2) at (1.175571,-1.618034) {};
\node [vertex] (v3) at (1.902113,0.618034) {};
\node [vertex] (v4) at (0.000000,2.000000) {};
\node [vertex] (v5) at (-0.951057,0.309017) {};
\node [vertex] (v6) at (-0.587785,-0.809017) {};
\node [vertex] (v7) at (0.587785,-0.809017) {};
\node [vertex] (v8) at (0.951057,0.309017) {};
\node [vertex] (v9) at (0.000000,1.000000) {};
\node [vertex] (v10) at (-0.769421,-0.250000) {};
\node [vertex] (v11) at (0.769421,-0.250000) {};
\node [vertex] (v12) at (0.000000,0.539345) {};
\node [left] at (v10) {\small $v_{12}$\,};
\draw (v0)--(v1);
\draw (v0)--(v4);
\draw (v0)--(v5);
\draw (v1)--(v2);
\draw (v1)--(v6);
\draw (v2)--(v3);
\draw (v2)--(v7);
\draw (v3)--(v4);
\draw (v3)--(v8);
\draw (v4)--(v9);
\draw (v5)--(v10);
\draw (v5)--(v12);
\draw (v6)--(v10);
\draw (v7)--(v11);
\draw (v8)--(v11);
\draw (v8)--(v12);
\draw (v9)--(v12);
\node [redvertex] at (v0) {};
\node [redvertex] at (v5) {};
\node [redvertex] at (v10) {};
\node [redvertex] at (v3) {};
\node [redvertex] at (v8) {};
\draw [ultra thick, red] (v0)--(v5)--(v10);
\draw [ultra thick, red] (v3)--(v8);
\node at (0,-2.25) {$X_{9}$};
\end{tikzpicture}
\hspace{1cm}
\begin{tikzpicture}[scale=1.1]
\tikzstyle{vertex}=[fill=blue!25, circle, draw=black,inner sep=0.7mm]
\tikzstyle{redvertex}=[fill=red, circle, draw=black,inner sep=0.7mm]
\node [vertex] (v0) at (-1.175571,-1.618034) {};
\node [vertex] (v1) at (1.175571,-1.618034) {};
\node [vertex] (v2) at (1.902113,0.618034) {};
\node [vertex] (v3) at (0.000000,2.000000) {};
\node [vertex] (v4) at (-1.902113,0.618034) {};
\node [vertex] (v5) at (-0.587785,-0.809017) {};
\node [vertex] (v6) at (0.587785,-0.809017) {};
\node [vertex] (v7) at (0.951057,0.309017) {};
\node [vertex] (v8) at (0.000000,1.000000) {};
\node [vertex] (v9) at (-0.951057,0.309017) {};
\node [vertex] (v10) at (0.000000,-0.809017) {};
\node [vertex] (v11) at (-0.769421,-0.250000) {};
\node [vertex] (v12) at (0.769421,-0.250000) {};
\node [vertex] (v13) at (0.000000,0.539345) {};
\node [left] at (v11) {\small $v_{12}$\,};
\draw (v0)--(v1);
\draw (v0)--(v4);
\draw (v0)--(v5);
\draw (v1)--(v2);
\draw (v1)--(v6);
\draw (v2)--(v3);
\draw (v2)--(v7);
\draw (v3)--(v4);
\draw (v3)--(v8);
\draw (v4)--(v9);
\draw (v5)--(v10);
\draw (v5)--(v11);
\draw (v6)--(v10);
\draw (v6)--(v12);
\draw (v7)--(v12);
\draw (v7)--(v13);
\draw (v8)--(v13);
\draw (v9)--(v11);
\draw (v9)--(v13);
\node [redvertex] at (v4) {};
\node [redvertex] at (v9) {};
\node [redvertex] at (v11) {};
\node [redvertex] at (v2) {};
\node [redvertex] at (v7) {};
\draw [ultra thick, red] (v2)--(v7);
\draw [ultra thick, red] (v4)--(v9)--(v11);
\node at (0,-2.25) {$X_{11}$};
\end{tikzpicture}
\end{center}
\caption{Ruling out $X_9$ and $X_{11}$}
\label{fig:x9x11}
\end{figure}

\begin{lemma}\label{notx9x11}
The subgraph $\widetilde{S} \not= X_9$, $X_{11}$.
\end{lemma}

\begin{proof}
In each case, let 
\[
T = \{u_1, v_1, v_{12}, u_4, v_4\}
\]
where $v_{12}$ is the common neighbour of $v_1$ and $v_2$.
Then 
\[
M_{TT}= 
\kbordermatrix{
& u_1 & v_1 & v_{12} & u_4 & v_4 \\
u_1 & 3 & 2 & 1 & 1 & 0\\
v_1 & 2 & 3 & 2 & 0 & 1\\
v_{12} & 1 & 2 & 3 & 0 & 0 \\
u_4 & 1 & 0 & 0 & 3 & 2\\
v_4 & 0 & 1 & 0 & 2 & 3
}
\]
which has determinant $-12$ and therefore $M$ is not positive semidefinite.
\end{proof}

We can now complete the proof of \cref{thm:girth5} which we restate here.

\begin{theorem*}(Restated \cref{thm:girth5})
If $G$ is a cubic graph of girth $5$ with no eigenvalues in $(-2,0)$ then $G$ is either the Petersen graph or the dodecahedron.
\end{theorem*}
\begin{proof}
If $C$ is a cycle of length $5$ in $G$ and 
$S = V(\cor{C})$, then by \cref{notx0x2x3x4x8x10x12,notx5x6,notx9x11}, the subgraph $\widetilde{S}$ is isomorphic to $X_1$ or $X_7$. In the first case, $\widetilde{S}$ is already cubic and is isomorphic to the Petersen graph.  

In the second case $\widetilde{S}$ is isomorphic to $X_7$ and this must hold true for \emph{every} choice of $C$ and its corresponding $S = V(\cor{C})$. Thus every $5$-cycle $C$ determines a corona with $5$ vertices of degree one that are each connected by a $2$-path to the previous/next vertex in the cyclic order determined by $C$. In addition, there are no other edges or 2-paths connecting the other pairs of vertices of degree one.

With $\widetilde{S}$ labelled as in \cref{fig:x7analysis} there are just five vertices of degree less than three, namely 
$V' = \{v_{01},v_{12},v_{23},v_{34},v_{04}\}$. In principle, there may be edges or 2-paths connecting vertices from $V'$, so first we rule these out. Without loss of generality we may assume that the endpoints of such an edge, or 2-path are $v_{01}$ and $v_{23}$. The two diagrams of \cref{fig:x7analysis} show the two cases with $v_{01}$ and $v_{23}$ adjacent in the first diagram and connected by a 2-path in the second. 
But in both cases, the corona of the $5$-cycle $(u_1, u_2, v_2, v_{12}, v_1, u_1)$ (pictured in red) does not lead to $X_7$ because $v_{01}$ and $v_{23}$ are degree one vertices in this corona that cannot be connected by either an edge or a 2-path. 

Therefore each of the five vertices in $V'$ is adjacent to another vertex and these vertices are distinct. We denote these vertices by $\{w_{01}, w_{12}, w_{23}, w_{34}, w_{04}\}$ where $w_{ij}$ is adjacent to $v_{ij}$. This results in the configuration shown in the first diagram of \cref{fig:xy}.

Now consider the corona of the highlighted $5$-cycle in \cref{fig:xy}. The vertices $w_{23}$ and $w_{12}$ are neighbouring vertices of degree one in this corona and so they must be joined by a $2$-path, which necessarily uses $w_{23}$ and hence $w_{12}$ and $w_{23}$ are adjacent in $G$. Repeating this argument for the corona of each of the five outer $5$-cycles yields the dodecahedron. 
\end{proof}

\begin{figure}
\begin{center}
\begin{tikzpicture}[scale=1.35]
\tikzstyle{vertex}=[fill=blue!25, circle, draw=black,inner sep=0.9mm]
\tikzstyle{redvertex}=[fill=red, circle, draw=black,inner sep=0.9mm]
\foreach \x in {0,2,4,6,8} {
\node [vertex] (v\x) at (90+36*\x:1.25cm) {};
}
\foreach \x in {1,3,5,7,9} {
\node [vertex] (v\x) at (90+36*\x:1.01cm) {};
}
\foreach \x in {0,1,...,4} {
\node [vertex] (u\x) at (90+72*\x:2.5cm) {};
}
\draw (v0)--(v1)--(v2)--(v3)--(v4)--(v5)--(v6)--(v7)--(v8)--(v9)--(v0);
\draw (u0)--(v0);
\draw (u1)--(v2);
\draw (u2)--(v4);
\draw (u3)--(v6);
\draw (u4)--(v8);
\draw (u0)--(u1)--(u2)--(u3)--(u4)--(u0);

\draw (v1)--(v5);

\node [redvertex] at (u1) {};
\node [redvertex] at (u2) {};
\node [redvertex] at (v4) {};
\node [redvertex] at (v3) {};
\node [redvertex] at (v2) {};

\draw [ultra thick, red] (u1)--(u2)--(v4)--(v3)--(v2)--(u1);

\node [above left] at (v1) {\small $v_{01}$};
\node [below left] at (v3) {\small $v_{12}$};
\node [below] at (v5.south) {\small $v_{23}$};
\node [below right] at (v7) {\small $v_{34}$};
\node [above right] at (v9) {\small $v_{04}$};
\end{tikzpicture}
\hspace{1cm}
\begin{tikzpicture}[scale=1.35]
\tikzstyle{vertex}=[fill=blue!25, circle, draw=black,inner sep=0.9mm]
\tikzstyle{redvertex}=[fill=red, circle, draw=black,inner sep=0.9mm]
\foreach \x in {0,2,4,6,8} {
\node [vertex] (v\x) at (90+36*\x:1.25cm) {};
}
\foreach \x in {1,3,5,7,9} {
\node [vertex] (v\x) at (90+36*\x:1.01cm) {};
}
\foreach \x in {0,1,...,4} {
\node [vertex] (u\x) at (90+72*\x:2.5cm) {};
}
\draw (v0)--(v1)--(v2)--(v3)--(v4)--(v5)--(v6)--(v7)--(v8)--(v9)--(v0);
\draw (u0)--(v0);
\draw (u1)--(v2);
\draw (u2)--(v4);
\draw (u3)--(v6);
\draw (u4)--(v8);
\draw (u0)--(u1)--(u2)--(u3)--(u4)--(u0);

\node [above left] at (v1) {\small $v_{01}$};
\node [below left] at (v3) {\small $v_{12}$};
\node [below] at (v5.south) {\small $v_{23}$};
\node [below right] at (v7) {\small $v_{34}$};
\node [above right] at (v9) {\small $v_{04}$};

\node [vertex] (x) at ($(v1)!0.5!(v5)$) {};
\draw (v1)--(x)--(v5);

\node [redvertex] at (u1) {};
\node [redvertex] at (u2) {};
\node [redvertex] at (v4) {};
\node [redvertex] at (v3) {};
\node [redvertex] at (v2) {};

\draw [ultra thick, red] (u1)--(u2)--(v4)--(v3)--(v2)--(u1);

%\draw [rotate = 18] (-1,-1) ellipse (1.5cm and 2cm);

\end{tikzpicture}
\end{center}
\caption{Cycles whose corona cannot be $X_7$.}
\label{fig:x7analysis}
\end{figure}

\begin{figure}
\begin{center}
\begin{tikzpicture}[scale=1.7]
\tikzstyle{vertex}=[fill=blue!25, circle, draw=black,inner sep=0.8mm]
\tikzstyle{redvertex}=[fill=red, circle, draw=black,inner sep=0.8mm]

\foreach \x in {0,2,4,6,8} {
\node [vertex] (v\x) at (90+36*\x:1.25cm) {};
}
\foreach \x in {1,3,5,7,9} {
\node [vertex] (v\x) at (90+36*\x:1.01cm) {};
}
\foreach \x in {0,1,...,4} {
\node [vertex] (u\x) at (90+72*\x:2.25cm) {};
}

\foreach \x in {1,3,5,7,9} {
\node [vertex] (x\x) at (90+36*\x:0.5cm) {};
\draw (x\x)--(v\x);
}

\node [right] at (x3) {\small $w_{12}$};
\node [right] at (x5) {\small $w_{23}$};
\node [below] at (v5.south) {\small $v_{23}$};

\draw (v0)--(v1)--(v2)--(v3)--(v4)--(v5)--(v6)--(v7)--(v8)--(v9)--(v0);
\draw (u0)--(v0);
\draw (u1)--(v2);
\draw (u2)--(v4);
\draw (u3)--(v6);
\draw (u4)--(v8);
\draw (u0)--(u1)--(u2)--(u3)--(u4)--(u0);
\node [redvertex] at (u1) {};
\node [redvertex] at (u2) {};
\node [redvertex] at (v4) {};
\node [redvertex] at (v3) {};
\node [redvertex] at (v2) {};
\draw [ultra thick, red] (u1)--(u2)--(v4)--(v3)--(v2)--(u1);
\end{tikzpicture}
\hspace{1cm}
\begin{tikzpicture}[scale=1.7]
\tikzstyle{vertex}=[fill=blue!25, circle, draw=black,inner sep=0.8mm]

\foreach \x in {0,2,4,6,8} {
\node [vertex] (v\x) at (90+36*\x:1.25cm) {};
}
\foreach \x in {1,3,5,7,9} {
\node [vertex] (v\x) at (90+36*\x:1.01cm) {};
}
\foreach \x in {0,1,...,4} {
\node [vertex] (u\x) at (90+72*\x:2.25cm) {};
}

\foreach \x in {1,3,5,7,9} {
\node [vertex] (x\x) at (90+36*\x:0.5cm) {};
\draw (x\x)--(v\x);
}

\draw (x1)--(x3)--(x5)--(x7)--(x9)--(x1);

\draw (v0)--(v1)--(v2)--(v3)--(v4)--(v5)--(v6)--(v7)--(v8)--(v9)--(v0);
\draw (u0)--(v0);
\draw (u1)--(v2);
\draw (u2)--(v4);
\draw (u3)--(v6);
\draw (u4)--(v8);
\draw (u0)--(u1)--(u2)--(u3)--(u4)--(u0);

\end{tikzpicture}
\end{center}
\caption{The dodecahedron is the only graph that arises}
\label{fig:xy}
\end{figure}

%% file: X0.tex
\begin{tikzpicture}[scale=0.8]
\tikzstyle{vertex}=[fill=blue!25, circle, draw=black,inner sep=0.7mm]
\node [vertex] (v0) at (0.000000,2.000000) {};
\node [vertex] (v1) at (-1.902113,0.618034) {};
\node [vertex] (v2) at (-1.175571,-1.618034) {};
\node [vertex] (v3) at (1.175571,-1.618034) {};
\node [vertex] (v4) at (1.902113,0.618034) {};
\node [vertex] (v5) at (0.000000,1.000000) {};
\node [vertex] (v6) at (-0.951057,0.309017) {};
\node [vertex] (v7) at (-0.587785,-0.809017) {};
\node [vertex] (v8) at (0.587785,-0.809017) {};
\node [vertex] (v9) at (0.951057,0.309017) {};
\draw (v0)--(v1);
\draw (v0)--(v4);
\draw (v0)--(v5);
\draw (v1)--(v2);
\draw (v1)--(v6);
\draw (v2)--(v3);
\draw (v2)--(v7);
\draw (v3)--(v4);
\draw (v3)--(v8);
\draw (v4)--(v9);
\node at (0,-2.25) {$X_{0}$};
\end{tikzpicture}

%% file: X1.tex
\begin{tikzpicture}[scale=0.8]
\tikzstyle{vertex}=[fill=blue!25, circle, draw=black,inner sep=0.7mm]
\node [vertex] (v0) at (0.000000,2.000000) {};
\node [vertex] (v1) at (-1.902113,0.618034) {};
\node [vertex] (v2) at (-1.175571,-1.618034) {};
\node [vertex] (v3) at (1.175571,-1.618034) {};
\node [vertex] (v4) at (1.902113,0.618034) {};
\node [vertex] (v5) at (0.000000,1.000000) {};
\node [vertex] (v6) at (-0.951057,0.309017) {};
\node [vertex] (v7) at (-0.587785,-0.809017) {};
\node [vertex] (v8) at (0.587785,-0.809017) {};
\node [vertex] (v9) at (0.951057,0.309017) {};
\draw (v0)--(v1);
\draw (v0)--(v4);
\draw (v0)--(v5);
\draw (v1)--(v2);
\draw (v1)--(v6);
\draw (v2)--(v3);
\draw (v2)--(v7);
\draw (v3)--(v4);
\draw (v3)--(v8);
\draw (v4)--(v9);
\draw (v5)--(v7);
\draw (v5)--(v8);
\draw (v6)--(v8);
\draw (v6)--(v9);
\draw (v7)--(v9);
\node at (0,-2.25) {$X_{1}$};
\end{tikzpicture}

%% file: X2.tex
\begin{tikzpicture}[scale=0.8]
\tikzstyle{vertex}=[fill=blue!25, circle, draw=black,inner sep=0.7mm]
\node [vertex] (v0) at (0.000000,2.000000) {};
\node [vertex] (v1) at (-1.902113,0.618034) {};
\node [vertex] (v2) at (-1.175571,-1.618034) {};
\node [vertex] (v3) at (1.175571,-1.618034) {};
\node [vertex] (v4) at (1.902113,0.618034) {};
\node [vertex] (v5) at (0.000000,1.000000) {};
\node [vertex] (v6) at (-0.951057,0.309017) {};
\node [vertex] (v7) at (-0.587785,-0.809017) {};
\node [vertex] (v8) at (0.587785,-0.809017) {};
\node [vertex] (v9) at (0.951057,0.309017) {};
\node [vertex] (v10) at (-0.475528,0.654508) {};
\draw (v0)--(v1);
\draw (v0)--(v4);
\draw (v0)--(v5);
\draw (v1)--(v2);
\draw (v1)--(v6);
\draw (v2)--(v3);
\draw (v2)--(v7);
\draw (v3)--(v4);
\draw (v3)--(v8);
\draw (v4)--(v9);
\draw (v5)--(v10);
\draw (v6)--(v10);
\node at (0,-2.25) {$X_{2}$};
\end{tikzpicture}

%% file: X3.tex
\begin{tikzpicture}[scale=0.8]
\tikzstyle{vertex}=[fill=blue!25, circle, draw=black,inner sep=0.7mm]
\node [vertex] (v0) at (0.000000,2.000000) {};
\node [vertex] (v1) at (-1.902113,0.618034) {};
\node [vertex] (v2) at (-1.175571,-1.618034) {};
\node [vertex] (v3) at (1.175571,-1.618034) {};
\node [vertex] (v4) at (1.902113,0.618034) {};
\node [vertex] (v5) at (0.000000,1.000000) {};
\node [vertex] (v6) at (-0.951057,0.309017) {};
\node [vertex] (v7) at (-0.587785,-0.809017) {};
\node [vertex] (v8) at (0.587785,-0.809017) {};
\node [vertex] (v9) at (0.951057,0.309017) {};
\node [vertex] (v10) at (-0.475528,0.654508) {};
\node [vertex] (v11) at (0.475528,0.654508) {};
\draw (v0)--(v1);
\draw (v0)--(v4);
\draw (v0)--(v5);
\draw (v1)--(v2);
\draw (v1)--(v6);
\draw (v2)--(v3);
\draw (v2)--(v7);
\draw (v3)--(v4);
\draw (v3)--(v8);
\draw (v4)--(v9);
\draw (v5)--(v10);
\draw (v5)--(v11);
\draw (v6)--(v10);
\draw (v9)--(v11);
\node at (0,-2.25) {$X_{3}$};
\end{tikzpicture}

%% file: X4.tex
\begin{tikzpicture}[scale=0.8]
\tikzstyle{vertex}=[fill=blue!25, circle, draw=black,inner sep=0.7mm]
\node [vertex] (v0) at (0.000000,2.000000) {};
\node [vertex] (v1) at (-1.902113,0.618034) {};
\node [vertex] (v2) at (-1.175571,-1.618034) {};
\node [vertex] (v3) at (1.175571,-1.618034) {};
\node [vertex] (v4) at (1.902113,0.618034) {};
\node [vertex] (v5) at (0.000000,1.000000) {};
\node [vertex] (v6) at (-0.951057,0.309017) {};
\node [vertex] (v7) at (-0.587785,-0.809017) {};
\node [vertex] (v8) at (0.587785,-0.809017) {};
\node [vertex] (v9) at (0.951057,0.309017) {};
\node [vertex] (v10) at (-0.475528,0.654508) {};
\node [vertex] (v11) at (0.475528,0.654508) {};
\node [vertex] (v12) at (-0.769421,-0.250000) {};
\draw (v0)--(v1);
\draw (v0)--(v4);
\draw (v0)--(v5);
\draw (v1)--(v2);
\draw (v1)--(v6);
\draw (v2)--(v3);
\draw (v2)--(v7);
\draw (v3)--(v4);
\draw (v3)--(v8);
\draw (v4)--(v9);
\draw (v5)--(v10);
\draw (v5)--(v11);
\draw (v6)--(v10);
\draw (v6)--(v12);
\draw (v7)--(v12);
\draw (v9)--(v11);
\node at (0,-2.25) {$X_{4}$};
\end{tikzpicture}

%% file: X5.tex
\begin{tikzpicture}[scale=0.8]
\tikzstyle{vertex}=[fill=blue!25, circle, draw=black,inner sep=0.7mm]
\node [vertex] (v0) at (0.000000,2.000000) {};
\node [vertex] (v1) at (-1.902113,0.618034) {};
\node [vertex] (v2) at (-1.175571,-1.618034) {};
\node [vertex] (v3) at (1.175571,-1.618034) {};
\node [vertex] (v4) at (1.902113,0.618034) {};
\node [vertex] (v5) at (0.000000,1.000000) {};
\node [vertex] (v6) at (-0.951057,0.309017) {};
\node [vertex] (v7) at (-0.587785,-0.809017) {};
\node [vertex] (v8) at (0.587785,-0.809017) {};
\node [vertex] (v9) at (0.951057,0.309017) {};
\node [vertex] (v10) at (-0.475528,0.654508) {};
\node [vertex] (v11) at (0.475528,0.654508) {};
\node [vertex] (v12) at (0.000000,-0.809017) {};
\draw (v0)--(v1);
\draw (v0)--(v4);
\draw (v0)--(v5);
\draw (v1)--(v2);
\draw (v1)--(v6);
\draw (v2)--(v3);
\draw (v2)--(v7);
\draw (v3)--(v4);
\draw (v3)--(v8);
\draw (v4)--(v9);
\draw (v5)--(v10);
\draw (v5)--(v11);
\draw (v6)--(v10);
\draw (v7)--(v12);
\draw (v8)--(v12);
\draw (v9)--(v11);
\node at (0,-2.25) {$X_{5}$};
\end{tikzpicture}

%% file: X6.tex
\begin{tikzpicture}[scale=0.8]
\tikzstyle{vertex}=[fill=blue!25, circle, draw=black,inner sep=0.7mm]
\node [vertex] (v0) at (0.000000,2.000000) {};
\node [vertex] (v1) at (-1.902113,0.618034) {};
\node [vertex] (v2) at (-1.175571,-1.618034) {};
\node [vertex] (v3) at (1.175571,-1.618034) {};
\node [vertex] (v4) at (1.902113,0.618034) {};
\node [vertex] (v5) at (0.000000,1.000000) {};
\node [vertex] (v6) at (-0.951057,0.309017) {};
\node [vertex] (v7) at (-0.587785,-0.809017) {};
\node [vertex] (v8) at (0.587785,-0.809017) {};
\node [vertex] (v9) at (0.951057,0.309017) {};
\node [vertex] (v10) at (-0.475528,0.654508) {};
\node [vertex] (v11) at (0.475528,0.654508) {};
\node [vertex] (v12) at (-0.769421,-0.250000) {};
\node [vertex] (v13) at (0.000000,-0.809017) {};
\draw (v0)--(v1);
\draw (v0)--(v4);
\draw (v0)--(v5);
\draw (v1)--(v2);
\draw (v1)--(v6);
\draw (v2)--(v3);
\draw (v2)--(v7);
\draw (v3)--(v4);
\draw (v3)--(v8);
\draw (v4)--(v9);
\draw (v5)--(v10);
\draw (v5)--(v11);
\draw (v6)--(v10);
\draw (v6)--(v12);
\draw (v7)--(v12);
\draw (v7)--(v13);
\draw (v8)--(v13);
\draw (v9)--(v11);
\node at (0,-2.25) {$X_{6}$};
\end{tikzpicture}

%% file: X7.tex
\begin{tikzpicture}[scale=0.8]
\tikzstyle{vertex}=[fill=blue!25, circle, draw=black,inner sep=0.7mm]
\node [vertex] (v0) at (0.000000,2.000000) {};
\node [vertex] (v1) at (-1.902113,0.618034) {};
\node [vertex] (v2) at (-1.175571,-1.618034) {};
\node [vertex] (v3) at (1.175571,-1.618034) {};
\node [vertex] (v4) at (1.902113,0.618034) {};
\node [vertex] (v5) at (0.000000,1.000000) {};
\node [vertex] (v6) at (-0.951057,0.309017) {};
\node [vertex] (v7) at (-0.587785,-0.809017) {};
\node [vertex] (v8) at (0.587785,-0.809017) {};
\node [vertex] (v9) at (0.951057,0.309017) {};
\node [vertex] (v10) at (-0.475528,0.654508) {};
\node [vertex] (v11) at (0.475528,0.654508) {};
\node [vertex] (v12) at (-0.769421,-0.250000) {};
\node [vertex] (v13) at (0.000000,-0.809017) {};
\node [vertex] (v14) at (0.769421,-0.250000) {};
\draw (v0)--(v1);
\draw (v0)--(v4);
\draw (v0)--(v5);
\draw (v1)--(v2);
\draw (v1)--(v6);
\draw (v2)--(v3);
\draw (v2)--(v7);
\draw (v3)--(v4);
\draw (v3)--(v8);
\draw (v4)--(v9);
\draw (v5)--(v10);
\draw (v5)--(v11);
\draw (v6)--(v10);
\draw (v6)--(v12);
\draw (v7)--(v12);
\draw (v7)--(v13);
\draw (v8)--(v13);
\draw (v8)--(v14);
\draw (v9)--(v11);
\draw (v9)--(v14);
\node at (0,-2.25) {$X_{7}$};
\end{tikzpicture}

%% file: X8.tex
\begin{tikzpicture}[scale=0.8]
\tikzstyle{vertex}=[fill=blue!25, circle, draw=black,inner sep=0.7mm]
\node [vertex] (v0) at (0.000000,2.000000) {};
\node [vertex] (v1) at (-1.902113,0.618034) {};
\node [vertex] (v2) at (-1.175571,-1.618034) {};
\node [vertex] (v3) at (1.175571,-1.618034) {};
\node [vertex] (v4) at (1.902113,0.618034) {};
\node [vertex] (v5) at (0.000000,1.000000) {};
\node [vertex] (v6) at (-0.951057,0.309017) {};
\node [vertex] (v7) at (-0.587785,-0.809017) {};
\node [vertex] (v8) at (0.587785,-0.809017) {};
\node [vertex] (v9) at (0.951057,0.309017) {};
\node [vertex] (v10) at (-0.475528,0.654508) {};
\node [vertex] (v11) at (0.293893,0.095492) {};
\node [vertex] (v12) at (0.000000,-0.809017) {};
\draw (v0)--(v1);
\draw (v0)--(v4);
\draw (v0)--(v5);
\draw (v1)--(v2);
\draw (v1)--(v6);
\draw (v2)--(v3);
\draw (v2)--(v7);
\draw (v3)--(v4);
\draw (v3)--(v8);
\draw (v4)--(v9);
\draw (v5)--(v10);
\draw (v5)--(v11);
\draw (v6)--(v10);
\draw (v7)--(v12);
\draw (v8)--(v11);
\draw (v8)--(v12);
\node at (0,-2.25) {$X_{8}$};
\end{tikzpicture}

%% file: X9rot.tex
\begin{tikzpicture}[scale=0.8]
\tikzstyle{vertex}=[fill=blue!25, circle, draw=black,inner sep=0.7mm]
\node [vertex] (v0) at (-1.902113,0.618034) {};
\node [vertex] (v1) at (-1.175571,-1.618034) {};
\node [vertex] (v2) at (1.175571,-1.618034) {};
\node [vertex] (v3) at (1.902113,0.618034) {};
\node [vertex] (v4) at (0.000000,2.000000) {};
\node [vertex] (v5) at (-0.951057,0.309017) {};
\node [vertex] (v6) at (-0.587785,-0.809017) {};
\node [vertex] (v7) at (0.587785,-0.809017) {};
\node [vertex] (v8) at (0.951057,0.309017) {};
\node [vertex] (v9) at (0.000000,1.000000) {};
\node [vertex] (v10) at (-0.769421,-0.250000) {};
\node [vertex] (v11) at (0.769421,-0.250000) {};
\node [vertex] (v12) at (0.000000,0.539345) {};
\draw (v0)--(v1);
\draw (v0)--(v4);
\draw (v0)--(v5);
\draw (v1)--(v2);
\draw (v1)--(v6);
\draw (v2)--(v3);
\draw (v2)--(v7);
\draw (v3)--(v4);
\draw (v3)--(v8);
\draw (v4)--(v9);
\draw (v5)--(v10);
\draw (v5)--(v12);
\draw (v6)--(v10);
\draw (v7)--(v11);
\draw (v8)--(v11);
\draw (v8)--(v12);
\draw (v9)--(v12);
\node at (0,-2.25) {$X_{9}$};
\end{tikzpicture}

%% file: X10.tex
\begin{tikzpicture}[scale=0.8]
\tikzstyle{vertex}=[fill=blue!25, circle, draw=black,inner sep=0.7mm]
\node [vertex] (v0) at (0.000000,2.000000) {};
\node [vertex] (v1) at (-1.902113,0.618034) {};
\node [vertex] (v2) at (-1.175571,-1.618034) {};
\node [vertex] (v3) at (1.175571,-1.618034) {};
\node [vertex] (v4) at (1.902113,0.618034) {};
\node [vertex] (v5) at (0.000000,1.000000) {};
\node [vertex] (v6) at (-0.951057,0.309017) {};
\node [vertex] (v7) at (-0.587785,-0.809017) {};
\node [vertex] (v8) at (0.587785,-0.809017) {};
\node [vertex] (v9) at (0.951057,0.309017) {};
\node [vertex] (v10) at (-0.475528,0.654508) {};
\node [vertex] (v11) at (0.293893,0.095492) {};
\node [vertex] (v12) at (-0.769421,-0.250000) {};
\node [vertex] (v13) at (0.000000,-0.809017) {};
\draw (v0)--(v1);
\draw (v0)--(v4);
\draw (v0)--(v5);
\draw (v1)--(v2);
\draw (v1)--(v6);
\draw (v2)--(v3);
\draw (v2)--(v7);
\draw (v3)--(v4);
\draw (v3)--(v8);
\draw (v4)--(v9);
\draw (v5)--(v10);
\draw (v5)--(v11);
\draw (v6)--(v10);
\draw (v6)--(v12);
\draw (v7)--(v12);
\draw (v7)--(v13);
\draw (v8)--(v11);
\draw (v8)--(v13);
\node at (0,-2.25) {$X_{10}$};
\end{tikzpicture}

%% file: X11rot.tex
\begin{tikzpicture}[scale=0.8]
\tikzstyle{vertex}=[fill=blue!25, circle, draw=black,inner sep=0.7mm]
\node [vertex] (v0) at (-1.175571,-1.618034) {};
\node [vertex] (v1) at (1.175571,-1.618034) {};
\node [vertex] (v2) at (1.902113,0.618034) {};
\node [vertex] (v3) at (0.000000,2.000000) {};
\node [vertex] (v4) at (-1.902113,0.618034) {};
\node [vertex] (v5) at (-0.587785,-0.809017) {};
\node [vertex] (v6) at (0.587785,-0.809017) {};
\node [vertex] (v7) at (0.951057,0.309017) {};
\node [vertex] (v8) at (0.000000,1.000000) {};
\node [vertex] (v9) at (-0.951057,0.309017) {};
\node [vertex] (v10) at (0.000000,-0.809017) {};
\node [vertex] (v11) at (-0.769421,-0.250000) {};
\node [vertex] (v12) at (0.769421,-0.250000) {};
\node [vertex] (v13) at (0.000000,0.539345) {};
\draw (v0)--(v1);
\draw (v0)--(v4);
\draw (v0)--(v5);
\draw (v1)--(v2);
\draw (v1)--(v6);
\draw (v2)--(v3);
\draw (v2)--(v7);
\draw (v3)--(v4);
\draw (v3)--(v8);
\draw (v4)--(v9);
\draw (v5)--(v10);
\draw (v5)--(v11);
\draw (v6)--(v10);
\draw (v6)--(v12);
\draw (v7)--(v12);
\draw (v7)--(v13);
\draw (v8)--(v13);
\draw (v9)--(v11);
\draw (v9)--(v13);
\node at (0,-2.25) {$X_{11}$};
\end{tikzpicture}

%% file: X12.tex
\begin{tikzpicture}[scale=0.8]
\tikzstyle{vertex}=[fill=blue!25, circle, draw=black,inner sep=0.7mm]
\node [vertex] (v0) at (0.000000,2.000000) {};
\node [vertex] (v1) at (-1.902113,0.618034) {};
\node [vertex] (v2) at (-1.175571,-1.618034) {};
\node [vertex] (v3) at (1.175571,-1.618034) {};
\node [vertex] (v4) at (1.902113,0.618034) {};
\node [vertex] (v5) at (0.000000,1.000000) {};
\node [vertex] (v6) at (-0.951057,0.309017) {};
\node [vertex] (v7) at (-0.587785,-0.809017) {};
\node [vertex] (v8) at (0.587785,-0.809017) {};
\node [vertex] (v9) at (0.951057,0.309017) {};
\node [vertex] (v10) at (-0.769421,-0.250000) {};
\draw (v0)--(v1);
\draw (v0)--(v4);
\draw (v0)--(v5);
\draw (v1)--(v2);
\draw (v1)--(v6);
\draw (v2)--(v3);
\draw (v2)--(v7);
\draw (v3)--(v4);
\draw (v3)--(v8);
\draw (v4)--(v9);
\draw (v5)--(v7);
\draw (v5)--(v8);
\draw (v6)--(v8);
\draw (v6)--(v10);
\draw (v7)--(v10);
\node at (0,-2.25) {$X_{12}$};
\end{tikzpicture}

%% file: girth9plus.tex
\subsection{Girth greater than 8}

In this subsection we dispose of the case where the girth is $9$ or more.

\begin{lemma}
If $G$ is a cubic graph of girth at least $9$, then $G$ has an eigenvalue in the interval $(-2,0)$. 
\end{lemma}

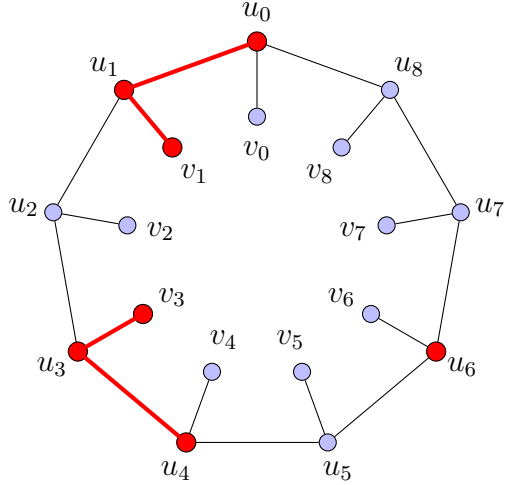
\begin{figure}
\begin{center}
\begin{tikzpicture}
\tikzstyle{vertex}=[fill=blue!25, circle, draw=black, inner sep=0.8mm]
\tikzstyle{redvertex}=[fill=red, circle, draw=black, inner sep=0.9mm]
\foreach \x in {0,1,...,8} {
\node [vertex](v\x) at (90+40*\x:1.75cm) {};
\node at (90+40*\x:1.3cm) {$v_{\x}$};
\node [vertex](u\x) at (90+40*\x:2.75cm) {};
\node at (90+40*\x:3.15cm) {$u_{\x}$};
\draw (v\x)--(u\x);
}

% \node [vertex](w11) at (90+36*11:2.75cm) {};
% \node at (90+36*11:2.25cm) {$w_{g-1}$};
% \node [vertex](v11) at (90+30*11:4.75cm) {};
% \node at (90+30*11:5.15cm) {$v_{g-1}$};
% \draw (v11)--(w11);

\draw (u0)--(u1)--(u2)--(u3)--(u4)--(u5)--(u6)--(u7)--(u8)--(u0);

\node [redvertex] at (v1) {};
\node [redvertex] at (u0) {};
\node [redvertex] at (u1) {};

\node [redvertex] at (v3) {};
\node [redvertex] at (u3) {};
\node [redvertex] at (u4) {};

\node [redvertex] at (u6) {};

%\draw [dashed] (v11) arc (60:-30:4.75);
%\draw [dashed] (w11) arc (60:-30:2.75);

\draw [ultra thick, red] (u0)--(u1)--(v1);
\draw [ultra thick, red] (v3)--(u3)--(u4);

\end{tikzpicture}
\end{center}
\caption{Configuration in the corona of a $9$-cycle}
\label{fig:corona9}
\end{figure}

\begin{proof}
Consider the corona of an $9$-cycle labelled as in \cref{fig:corona9}. As $G$ has girth $9$, the vertices $V = \{v_0, v_1, \ldots, v_8\}$ are distinct, pairwise non-adjacent, and no two have a common neighbour.

If we take $T = \{u_0,u_1,v_1,u_3,u_4,v_3,u_6\}$ as illustrated in red in \cref{fig:corona9} then
\[
M_{TT} =
\kbordermatrix{
&u_0&u_1&v_1&u_3&u_4&v_3&u_6\\
u_0&3&2&1&0&0&0&0\\
u_1&2&3&2&1&0&0&0\\
v_1&1&2&3&0&0&0&0\\
u_3&0&1&0&3&2&2&0\\
u_4&0&0&0&2&3&1&1\\
v_3&0&0&0&2&1&3&0\\
u_6&0&0&0&0&1&0&3\\
}
\]

% \[
% M_{TT} =
% \left[
% \begin{array}{ccc|ccc|c}
% 3&2&1&0&0&0&0\\
% 2&3&2&1&0&0&0\\
% 1&2&3&0&0&0&0\\
% \hline
% 0&1&0&3&2&2&0\\
% 0&0&0&2&3&1&1\\
% 0&0&0&2&1&3&0\\
% \hline
% 0&0&0&0&1&0&3\\
% \end{array}
% \right]
% \]
and because $\det M_{TT} = -16$ it follows that $M$ is not positive semidefinite and therefore $G$ has an eigenvalue in $(-2,0)$.

If $G$ has girth greater than $9$, then simply taking the same subset of vertices from the corona of any girth-cycle of $G$ will yield the same matrix $M_{TT}$.
\end{proof}

%% file: conclusion.tex
\section{Conclusion}

In studying spectral gap sets for cubic graphs, Koll\'ar and Sarnak \cite{KolSar2021} identified four maximal spectral gap intervals for cubic graphs. Two of these intervals, namely $(-1,1)$ and $(-2,0)$, have length $2$ which is the maximum possible length for a spectral gap interval for cubic graphs.

In \cite{guo2024cubicgraphseigenvaluesinterval} we used the fact that a cubic graph $G$ has no eigenvalues in $(-1,1)$ if and only if the matrix $(A(G)-I)(A(G)+I)$ is positive semidefinite to provide a complete characterisation of these graphs. This complete characterisation allowed us to show that $(-1,1)$ is a maximal gap \emph{set} thereby resolving a question left open in Koll\'ar and Sarnak \cite{KolSar2021}. 

% In the current paper, we have followed the same basic idea of associating a matrix to a cubic graph in such a way that the matrix is positive semidefinite if and only if the graph has no eigenvalues in $(-2,0)$, thereby completing the characterist of these two fundamental spectral gap intervals.

We note that the characterisation in
\cite{guo2024cubicgraphseigenvaluesinterval} is easier and more elegant than the one in this paper due to two serendipitous features related to the particular interval $(-1,1)$. Firstly, the fact that $(-1,1)$ is symmetric about the origin allowed us to reduce the problem to the bipartite case and secondly, a fortuitous numerical coincidence meant that we could directly use the well-known classification of graphs with minimum eigenvalue at least $-2$ to resolve one case.

It is conceivable that the same techniques could be used for other spectral gap intervals but, at least for cubic graphs, there are no other really natural spectral gap intervals left to analyse. In particular, a spectral gap interval for cubic graphs has length at most $2$ and the intervals $(-1,1)$ and $(-2,0)$ are only ones known of this length. There are no other spectral gap intervals of length $2$ with integer endpoints.

\subsection{Computational Note}

We have checked the results computationally for the cubic graphs on up to 26 vertices (there are $2094480864$ or just over 2 billion graphs of this order). We used Pari/GP \cite{PARI2}, which has various routines for exact computation with integer polynomials thereby avoiding even the remote possibility that floating-point numbers very very close to $0$ or $-2$ are incorrectly classified due to round-off error.

\section{Acknowledgements}

We thank Alicia Koll\'ar and Peter Sarnak both for their interest in these results, and for their assistance in explaining and clarifying various details of their work on spectral gap sets. 